\documentclass[a4paper,fleqn,12pt]{article}

\usepackage{vmargin}
\setpapersize{A4}
\setmarginsrb{35mm}{30mm}{30mm}{30mm}{0mm}{0mm}{10mm}{3\baselineskip}

\usepackage{graphicx}

\usepackage{amsmath}
\usepackage{amssymb}
\usepackage{theorem}

\newtheorem{satz}{Satz}[section]                    

\newtheorem{lemma}[satz]{Lemma}

\newtheorem{corollary}[satz]{Corollary}
\newtheorem{theorem}[satz]{Theorem}

\newtheorem{proposition}[satz]{Proposition}

{\theorembodyfont{\rmfamily} \newtheorem{remark}[satz]{Remark}}

\newenvironment{myen}[1]
	       {%
		 \begin{trivlist}%
		   \settowidth{\labelwidth}{\textit{#1:}}%
		   \setlength{\labelsep}{0in}%
		   \setlength{\itemindent}{\labelwidth}%
		 \item[\textit{#1:}]\ %
	       }
               {\end{trivlist}}
	       
\newenvironment{proof}{\begin{myen}{Proof}}{\end{myen}}
\newenvironment{myproof}{\begin{proof}}{\fillbox\end{proof}}
\newcommand{\fillbox}{\hspace*{\fill}\ensuremath{\Box}}

\numberwithin{equation}{section}

\makeatletter
\renewcommand{\section}{\@startsection%
  {section}%
  {1}%
  {0mm}%
  {-\baselineskip}%
  {0.5\baselineskip}%
  {\normalfont\large\bf}}%
\renewcommand{\subsection}{\@startsection%
  {subsection}%
  {2}%
  {0mm}%
  {-0.5\baselineskip}%
  {0.2\baselineskip}%
  {\normalfont\bf}}%
\makeatother

\newcounter{mynumcounter}
\newenvironment{mynum}%
	       {%
		 \begin{list}{(\roman{mynumcounter})\hspace*{\fill}}%
		   {
		     \setlength{\topsep}{0cm}
		     \setlength{\partopsep}{0cm}
		     \setlength{\itemsep}{0cm}
		     \setlength{\parsep}{0cm}
		     \setlength{\leftmargin}{0cm}
		     \setlength{\itemindent}{8mm}		   
		     \setlength{\labelsep}{3mm}
		     \setlength{\labelwidth}{5mm}
		     \usecounter{mynumcounter}
		   }%
	       }
	       {\end{list}}

\newcommand{\eps}{\varepsilon}
\newcommand{\del}{\partial}

\newcommand{\dd}[2]{\frac{\del #1}{\del #2}}
\newcommand{\ddeval}[3]{\left.\dd{#1}{#2}\right|_{#3}}

\newcommand{\DD}[2]{\frac{d #1}{d #2}}
\newcommand{\DDeval}[3]{\left.\DD{#1}{#2}\right|_{#3}}

\newcommand{\Ric}{\mathrm{Ric}}
\newcommand{\Scal}{\mathrm{Scal}}
\newcommand{\graph}{\mathrm{graph}}

\newcommand{\tr}{\mathrm{tr}}
\renewcommand{\div}{\mathrm{div}}

\newcommand{\IR}{\mathbf{R}}

\newcommand{\CF}{\mathcal{F}}
\newcommand{\CH}{\mathcal{H}}
\newcommand{\CL}{\mathcal{L}}

\newcommand{\CP}{\mathcal{P}}
\newcommand{\CR}{\mathcal{R}}
\newcommand{\CS}{\mathcal{S}}

\newcommand{\Id}{\mathrm{Id}}
\newcommand{\id}{\mathrm{id}}

\newcommand{\rmd}{\mathrm{d}}
\newcommand{\Acirc}{\hspace*{4pt}\raisebox{8.5pt}{\makebox[-4pt][l]{$\scriptstyle\circ$}}A}
\newcommand{\rmin}{r_\mathrm{min}}
\newcommand{\const}{\mathrm{const}}
\newcommand{\supp}{\mathrm{supp}}
\newcommand{\dist}{\mathrm{dist}}
\newcommand{\tensor}{\otimes}
\newcommand{\LHP}{L^{\CHP}}
\newcommand{\dmu}{\,\rmd\mu}

\newcommand{\ra}{\rangle}
\newcommand{\la}{\langle}

\newcommand{\gasym}{\|g-g^S\|_{C^2_{-1-\delta}}}

\newcommand{\Kasym}{\|K\|_{C^1_{-2-\delta}}}

\newcommand{\CHP}{\ensuremath{\mathcal{H}\pm\mathcal{P}}}
\newcommand{\Oeta}{O(\eta)}
\newcommand{\IAK}{\int_{A(k)}\!\!\!}

\begin{document}
  \begin{center}
    {\Large\bf Foliations of asymptotically flat 3-manifolds by 2-surfaces of
      prescribed mean curvature}\\[4mm]
    {\normalsize Jan Metzger\footnote{E-Mail:
      metzger@aei.mpg.de. Partially supported by the DFG,
      SFB 382 at T{\"u}bingen University.}}\\[2mm]
    { \small\it Max-Planck-Institut f\"ur Gravitationsphysik\\ 
      Am M\"uhlenberg 1, D-14476 Potsdam, Germany.}\\[2mm]
      {\normalsize}
  \end{center}
  \begin{quote}
    {\small{\bf Abstract.}  We construct 2-surfaces of prescribed mean
      curvature in 3-mani\-folds carrying asymptotically flat initial
      data for an isolated gravitating system with rather general
      decay conditions. The surfaces in question form a regular
      foliation of the asymptotic region of such a manifold. We
      recover physically relevant data, especially the ADM-momentum,
      from the geometry of the foliation.
      
      For a given set of data $(M,g,K)$, with a three dimensional
      manifold $M$, its Riemannian metric $g$, and the second
      fundamental form $K$ in the surrounding four dimensional Lorentz
      space time manifold, the equation we solve is $H+P=\const$ or
      $H-P=\const$. Here $H$ is the mean curvature, and $P = tr K$ is
      the 2-trace of $K$ along the solution surface.  This is a
      degenerate elliptic equation for the position of the surface. It
      prescribes the mean curvature anisotropically, since $P$ depends
      on the direction of the normal.}
  \end{quote}
  \section{Introduction and Statement of Results}
Surfaces with prescribed mean curvature play an important role for
example in the field of general relativity. Slicings are frequently
used to find canonic objects simplifying the treatment of the four
dimensional space-time. A prominent setting is the ADM 3+1
decomposition \cite{ADM:1961} of a four dimensional manifold into
three dimensional spacelike slices. Such slices are often chosen by
prescribing their mean curvature in the four geometry. In contrast, we
consider the subsequent slicing of a three dimensional spacelike slice
by two dimensional spheres with prescribed mean curvature in the three
geometry.

To be more precise, let $(M,g,K)$ be a set of initial data. That is,
$(M,g)$ is a three dimensional Riemannian manifold and $K$ is a
symmetric bilinear form on $M$. This can be interpreted as the
extrinsic curvature of $M$ in the surrounding four dimensional space
time. We consider 2-surfaces $\Sigma$ satisfying one of the
quasilinear degenerate elliptic equations $H\pm P=\const$ where $H$
ist the mean curvature of $\Sigma$ in $(M,g)$ and $P=\tr^\Sigma K$ is
the two dimensional trace of $K$.

In the case where $K\equiv 0$ this equation particularizes to
$H=\const$, which is the Euler-Lagrange equation of the isoperimetric
problem. This means that surfaces satisfying $H=\const$ are stationary
points of the area functional with respect to volume preserving
variations. Yau suggested to use such surfaces to describe physical
information in terms of geometrically defined objects. Indeed Huisken
and Yau \cite{HY:1996} have shown that the asymptotic end of an
asymptotically flat manifold, with appropriate decay conditions on the
metric, is uniquely foliated by such surfaces which are stable with
respect to the isoperimetric problem. The Hawking mass
\[ m_\text{H}(\Sigma) = \frac{|\Sigma|^{1/2}}{(16\pi)^{3/2}}\left(16\pi-
\int_\Sigma H^2\dmu \right) \] of such a surface $\Sigma$ is monotone
on this foliation and converges to the ADM-mass. This foliation can
also be used to define the center of mass of an isolated system since
for growing radius, the surfaces approach Euclidean spheres with a
converging center. Therefore the static physics of an isolated system
considered as point mass is contained in the geometry of the
$H=\const$ foliation. However, these surfaces are defined independently
of $K$, such that no dynamical physics can be found in their
geometry. A different proof of the existence of CMC surfaces is due to
Ye \cite{Ye:1996}.

The goal of this paper is to generalize the CMC foliations to include
the dynamical information into the definition of the foliation. The
equation $H\pm P=\const$ was chosen since apparent horizons satisfying
$H=0$ in the case $K\equiv 0$ generalize to surfaces satisfying $H\pm
P=0$ when $K$ does not necessarily vanish. We made this choice with
the Penrose inequality \cite{Penrose:1973} in mind. This inequality 
estimates the ADM-mass of an isolated system by the area of a black
hole horizon $\Sigma$
\[m_\text{ADM}\geq \sqrt{\frac{|\Sigma|}{16\pi}}\,.\]
In the case $K\equiv 0$ this becomes the Riemannian Penrose
inequality, which says that if $\Sigma$ is an outermost minimal
surface then the above inequality is valid. It was proved by Huisken
and Ilmanen \cite{HI:1997} and Bray \cite{Bray:1997}, both using
prescribed mean curvature surfaces. The proof of Bray generalizes the
situation in which an outermost minimal surface is part of the stable
CMC foliation from \cite{HY:1996}, in that it shows that the Hawking
mass is monotone on isoperimetric surfaces when their enclosed volume
and area increase even though they may not form a foliation. While a
fully general apparent horizon Penrose inequality does not seem to be
true \cite{Ben-Dov:2004}, generalizing this picture is of interest as
it may help to investigate whether there is a replacement which is
still true.

We consider asymptotically flat data describing isolated gravitating
systems. For constants $m>0$, $\delta\geq 0$, $\sigma\geq0$, and
$\eta\geq 0$ data $(M,g,K)$ will be called
$(m,\delta,\sigma,\eta)$-asymptotically flat if there exists a compact
set $B\subset M$ and a diffeomorphism $x: M\setminus B\to
\IR^3\setminus B_\sigma(0)$ such that in these coordinates $g$ is
asymptotic to the conformally flat spatial Schwarzschild metric $g^S$
representing a static black hole of mass $m$. Here, $g^S=\phi^4 g^e$,
where $\phi=1+\frac{m}{2r}$, $g^e$ is the Euclidean metric, and $r$ is
the Euclidean radius. The asymptotics we require for $g$ and $K$ are
\begin{gather}
  \label{eq:Einl-asympt-1}
  \sup_{\IR^3\setminus B_\sigma(0)} \left( r^{1+\delta} |g-g^S| + r^{2+\delta} |\nabla^g - \nabla^S| +
  r^{3+\delta}|\Ric^g -\Ric^S| \right) < \eta \,,\\
  \label{eq:Einl-asympt-2}  
  \sup_{\IR^3\setminus B_\sigma(0)} \left( r^{2+\delta} |K| +
  r^{3+\delta}|\nabla^g K| \right) < \eta\,.
\end{gather}
Here $\nabla^g$ and $\nabla^S$ denote the Levi-Civita connections of
$g$ and $g^S$ on $TM$, such that $\nabla^g-\nabla^S$ is a
$(1,2)$-tensor. Furthermore $\Ric^g$ and $\Ric^S$ denote the 
respective Ricci tensors of $g$ and $g^S$. That is, we consider data arising
from a perturbation of the Schwarzschild data $(g^S,0)$.

The main theorem will be proved for $\delta=0$ and $\eta=\eta(m)$
small compared to $m>0$. These conditions are optimal in the sense that
we only impose conditions on geometric quantities, not on partial
derivatives. They include far more general data than similar
results. Huisken and Yau \cite{HY:1996} for example demand that
$g-g^S$ decays like $r^{-2}$ with corresponding conditions on the
decay of the derivatives up to fourth order, while we only need
derivatives up to second order.  Christodoulou and Klainerman
\cite{CK:1993} use asymptotics with $g-g^S$ decaying like $r^{-3/2}$
with decay conditions on the derivatives up to fourth order, and $K$
like $r^{-5/2}$ with decay conditions on derivatives up to third
order, whereas our result needs to two levels of differentiability
less. In addition we allow data with nonzero ADM-momentum.
For such data with $\delta=0$ we can prove the following:
\begin{theorem}
  \label{thm:haupt}
  Given $m>0$ there is $\eta_0=\eta_0(m)>0$, such that if the data
  $(M,g,K)$ are $(m,0,\sigma,\eta_0)$-asymptotically flat for some
  $\sigma>0$, there is $h_0=h_0(m,\sigma)$ and a differentiable map
  \[ F: (0,h_0) \times S^2 \rightarrow M : (h,p) \mapsto F(h,p) \]
  satisfying the following statements.
  \begin{mynum}
  \item The map $F(h,\cdot):S^2\to M$ is an embedding. The surface
    $\Sigma_h = F(h,S^2)$ satisfies $H+P = h$ with respect to
    $(g,K)$. Each $\Sigma_h$ is convex $|A|^2\leq 4\det A$.
  \item There is a connected compact set $\bar B\subset M$, such that $F((0,h_0) ,
    S^2)=M\setminus \bar B$.
  \item The surfaces $F(h,S^2)$ form a regular foliation.
  \item There is a constant $C$ such that for all $h$ the surfaces $\Sigma_h$ satisfy
    \[ \|\nabla \Acirc\|_{L^2(\Sigma)}^2 + |\Sigma_h|^{-1}\|\Acirc\|_{L^2(\Sigma)}^2
    \leq C\eta^2|\Sigma|^{-2}\,. \]
  \item There are $\sup$-estimates for all curvature quantities on
     $\Sigma_h$, cf. section~\ref{s:apriori2}.
  \item Every convex surface $\Sigma$ with $H+P=h$ contained in
  $\IR^3\setminus B_{h^{-2/3}}(0)$ equals $\Sigma_h$. Hence the
  foliation is unique in the class of convex foliations.
  \end{mynum}
\end{theorem}
An analogous theorem holds for foliations with $H-P=\const$.

This theorem does not need that $(M,g,K)$ satisfy the constraint
equations. It can be generalized to give the existence of a foliation
satisfying $H+P_0(\nu)=\const$, where $P_0:SM\to\IR^3$ is a function
on the sphere bundle of $M$ with the same decay as $K$.

Our result includes the existence results from Huisken and Yau
\cite{HY:1996} for CMC foliations. Their uniqueness result for
individual surfaces can be proved in a smaller class, while the global
uniqueness result holds in the general case (cf. remark
\ref{rem:conditions}).

By rescaling $(g,K)$, the dependence of $\eta_0$ and $h_0$
on $m$ can be exposed. Consider the map $F_\sigma: x \mapsto \sigma x$,
and let $g_\sigma := \sigma^{-2}F_\sigma^* g^S$ and $K_\sigma:=
\sigma^{-1} F^*_\sigma K$. If $(g,K)$ is
$(m,0,\sigma,\eta)$-asymptotically flat, then $(g_{m},K_{m})$ is
$(1,0,m\sigma, m^{-1}\eta)$-asymptotically flat. Therefore
$\eta_0(m) = m \eta_0(1)$, and $h_0(\sigma,m) = m h_0(m\sigma,1)$.

Section~\ref{s:prelim} introduces some notation. In
sections~\ref{s:apriori1} and~\ref{s:apriori2} we carry out the a
priori estimates for the geometric quantities first in Sobolev norms
and then in the $\sup$-norm. Using these estimates we examine the
linearization of the operator $H\pm P$ in section~\ref{s:linearized}
and show that it is invertible. This is used in
section~\ref{s:foliation} to prove theorem~\ref{thm:haupt}. Finally,
in section~\ref{s:special} we use special asymptotics of
$(g,K)$ to investigate the connection between the foliation and the
linear momentum of these data.

  \section{Preliminaries}
\label{s:prelim}
\subsection{Notation}
Let $M$ be a three dimensional manifold. We will denote a
Riemannian metric on $M$ by $g$, or in coordinates by $g_{ij}$. Its inverse
is written as $g^{-1}=\{g^{ij}\}$. The Levi-Civita connection of $g$ is
denoted by $\nabla$, the Riemannian curvature tensor by $R$, the Ricci
tensor by $\Ric$, and the scalar curvature by $\Scal$.

Let $\Sigma$ be a hypersurface in $M$. Let $\gamma^g$ denote the
metric on $\Sigma$ induced by $g$, and let $\nu^g$ denote its normal. The
second fundamental form of $\Sigma$ is denoted by $A^g$, its mean
curvature by $H^g$, and the traceless part of $A^g$ by
$\Acirc^g=A^g-\frac{1}{2}H^g\gamma^g$.

We follow the Einstein summation convention and sum over Latin indices
from 1 to 3 and over Greek indices from 1 to 2.

We use the usual function spaces on compact surfaces with their usual
norms. The $L^p$-norm of an $(s,t)$-tensor $T$ with respect to the
metric $\gamma$ on $\Sigma$ is denoted by
\[ \|T\|^p_{L_{(s,t)}^p(\Sigma,\gamma)} = \int_\Sigma |T|_\gamma^p\dmu^\gamma\,. \]
The space $L_{(s,t)}^p(\Sigma)$ of $(s,t)$-tensors is the completion
of the space of smooth $(s,t)$-tensors with respect to this norm. In
the sequel we will drop the subscripts $(s,t)$, since norms will be
used unambiguously. The Sobolev-norm $W^{k,p}(\Sigma)$ is defined as
\[ \|T\|^p_{W^{k,p}(\Sigma)} = \|T\|^p_{L^p(\Sigma)}+\ldots+\|\nabla^k
T\|^p_{L^p(\Sigma)}\,, \] 
where $\nabla^K T$ is the $k$-th covariant derivative of $T$. Again,
the space $W^{k,p}(\Sigma)$ is the completion of the smooth tensors
with respect to this norm.

For a smooth tensor $T$, define the H\"older semi-norm  by 
\[ [T]_{p,\alpha} : = \sup_{p\neq q} \frac{ |P_q T(q) - T(p)|
}{\dist(p,q)^\alpha}\,, \]
where $P_q$ denotes parallel translation along the shortest geodesic
from $p$ to $q$, and the supremum is taken over all $p\neq q$ with
$\dist(p,q)$ less than the injectivity radius of $\Sigma$. Define the
H\"older norm $\|T\|_{C^{k,\alpha}(\Sigma)}$ as 
\[ \|T\|_{C^{k,\alpha}_{(s,t)}(\Sigma)} := \sup_\Sigma |T| + \sup_\Sigma
|\nabla T| + \cdots + \sup_\Sigma|\nabla^k T| + \sup_{p\in\Sigma} [\nabla^k
T]_{p,\alpha}\,. \] 
We assume in the following that $(M,g,K)$ and all hypersurfaces are
smooth, i.e. $C^\infty$. However, to prove theorem \ref{thm:haupt} it
is obviously enough to assume $g$ to be $C^2$ and $K$ to be $C^1$. The
a priori estimates from sections \ref{s:apriori1} and \ref{s:apriori2}
are valid for surfaces of class $W^{3,p}$, when $p$ is large enough.

\subsection{Extrinsic Geometry}
Let $\Sigma\subset(M,g)$ be a hypersurface. The second fundamental
form $A_{\alpha\beta}$ and the Riemannian curvature tensor
$\CR_{\alpha\beta\gamma\delta}$ of $\Sigma$ are connected to the curvature
$R_{ijkl}$ of $M$ via the Gauss and Codazzi equations
\begin{align}
  \label{eq:gauss}
  \CR_{\alpha\beta\delta\eps} 
  &= R_{\alpha\beta\delta\eps} +
  A_{\alpha\delta}A_{\beta\eps} - A_{\alpha\eps}A_{\beta\delta} 
  \\
  \label{eq:codazzi}
  R_{i\alpha\beta\delta}\nu^i 
  &= \nabla_\delta A_{\alpha\beta} -
  \nabla_\beta A_{\alpha\delta} \,.
\end{align}
Together, these imply the Simons identity \cite{Simons:1968,SSY:1975}
\begin{equation}   
  \label{eq:simons}
  \begin{split}
    \Delta^\Sigma A_{\alpha\beta} 
    &=     
    \nabla^\Sigma_\alpha \nabla^\Sigma_\beta H 
    + H A_\alpha^\delta A_{\delta\beta} 
    - |A|^2  A_{\alpha\beta} + A_\alpha^\delta R_{\eps\beta\eps\delta}
    + A^{\delta\eps}R_{\delta\alpha\beta\eps} 
    \\    
    & \phantom{=}
    + \nabla^\Sigma_\beta(\Ric_{\alpha k}\nu^k) + \nabla^{\Sigma\,\delta}
    (R_{k\alpha\beta\delta}\nu^k)  \,.
  \end{split}
\end{equation}
Note that the last two terms were not differentiated with the Leibniz
rule. Equation \eqref{eq:simons} therefore differs slightly from how
the Simons identity is usually stated.
\subsection{Round surfaces in Euclidean space}
The key tool in obtaining a priori estimates for the surfaces in
question is the following theorem by DeLellis and M\"uller
\cite[Theorem 1]{DLM:2003}. 
\begin{theorem}
  \label{thm:lm}
  There exists a universal constant $C$ such that for each compact
  connected surface without boundary $\Sigma$, with area
  $|\Sigma|=4\pi$, the following estimate holds:
  \begin{equation*}
    \| A - \gamma \|_{L^2(\Sigma)} \leq C \|\Acirc\|_{L^2(\Sigma)}\,.
  \end{equation*}
  If in addition $\|\Acirc\|_{L^2(\Sigma)}\leq 8\pi$, then $\Sigma$
  is a sphere, and there exists a conformal map
  $\psi:S^2\to\Sigma\subset \IR^3$ such that
  \begin{equation*}
    \| \psi - (a + \id_{S^2}) \|_{W^{2,2}(S^2)} \leq
    C\|\Acirc\|_{L^2(\Sigma)}\,,
  \end{equation*}
  where $\id_{S^2}$ is the standard embedding of $S^2$ onto the
  sphere $S_1(0)$ in $\IR^3$, and
  $a=|\Sigma|^{-1}\int_\Sigma\id_\Sigma\dmu$ is the center of
  gravity of $\Sigma$.
\end{theorem}
DeLellis and M\"uller \cite[3,6.1,6.3]{DLM:2003} also prove the
following useful estimates for the normal $\nu$ and the conformal
factor $h^2$ of such surfaces:
\begin{gather*}
  C^{-1} \leq h\leq C\,,\\
  \| h-1\|_{W^{1,2}(S^2)} \leq C\|\Acirc\|_{L^2(\Sigma)}\,, \\
  \| N - \nu\circ\psi \|_{W^{1,2}(S^2)}\leq C\|\Acirc\|_{L^2(\Sigma)}\,.
\end{gather*}
Here, $N$ is the normal of $S_1(a)$, and $h^2$ is the conformal factor
of $\psi$, such that if $\gamma^S$ denotes the metric on $S_R(a)$ and
$\gamma$ the metric on $\Sigma$, then we have $\psi^*\gamma =
h^2\gamma^S$.

To translate these inequalities into a scale invariant form for
surfaces which do not necessarily have area $|\Sigma|=4\pi$, we
introduce the Euclidean area radius
$R_e=\sqrt{|\Sigma|/4\pi}$. The first part of theorem
\ref{thm:lm} implies that for a general surface $\Sigma$
\begin{equation}
  \label{eq:lm_Aest}
  \| A - R_e^{-1}\gamma \|_{L^2(\Sigma)} \leq
  C\|\Acirc\|_{L^2(\Sigma)}\,.
\end{equation}
In the case $\|\Acirc\|_{L^2(\Sigma)}\leq 8\pi$, the second part of
theorem \ref{thm:lm} gives that there exist
$a_e:=|\Sigma|^{-1}\int_\Sigma\id_\Sigma\dmu\in\IR^3$ and a conformal
parameterization $\psi:S_{R_e}(a_e)\to\Sigma$. By the Sobolev
embedding on $S^2$ we obtain the following estimates for $\psi$, its
conformal factor $h^2$, and the difference of the normal $N$ of
$S_{R_e}(a_e)$ and the normal $\nu$ of $\Sigma$:
\begin{gather}
  \label{eq:lm_sup_est}
  {\textstyle \sup_{S_{R_e}(a_e)}}\left|\psi - \id_{S_{R_e}(a_e)}\right| \leq
  CR_e\|\Acirc\|_{L^2(\Sigma)}\,, \\
  \label{eq:lm_conf_est}
  \|h^2-1\|_{L^2(S_{R_e}(a_e))} \leq CR_e\|\Acirc\|_{L^2(\Sigma)}\,,\\
  \label{eq:lm_nu_est}
  \|N \circ \id_{S_{R_e}(a_e)} - \nu \circ\psi
  \|_{L^2(S_{R_e}(a_e))} \leq CR_e\|\Acirc\|_{L^2(\Sigma)}\,.
\end{gather}
\subsection{Asymptotically flat metrics}
Let $g^S$ be the spatial, conformally flat Schwarzschild metric on
$\IR^3\setminus\{0\}$. Namely, let $g^S_{ij} = \phi^4 g^e_{ij}$
with $\phi= 1 + \frac{m}{2r}$, where $g^e_{ij}=\delta_{ij}$ is the
Euclidean metric, and $r$ the Euclidean radius on $\IR^3$. Here and
in the sequel we will suppress the dependence of $g^S$ on the mass
parameter $m$. However, we will restrict ourselves to the case
$m>0$. The Ricci curvature of $g^S$ is given by
\begin{equation}
  \label{eq:schw_ric}
  \Ric_{ij}^S =
  \frac{m}{r^3}\phi^{-2}\left(\delta_{ij}-3\rho_i\rho_j\right)\,,
\end{equation}
where $\rho=x/r$ is the radial vector field on
$\IR^3$, whence $\Scal^S=0$.

Omiting $K$ and saying that data $(M,g)$ are
$(m,\delta,\sigma,\eta)$-asymptotically flat, we mean that $K\equiv 0$
and $(M,g,K\equiv 0)$ is $(m,\delta,\sigma,\eta)$-asymptotically
flat. Recall that then there exists a compact set $B\subset M$ and a
diffeomorphism $x:M\setminus B\to \IR^3\setminus B_\sigma(0)$, such
that in these coordinates the following 'norm'
\begin{equation}
  \label{eq:gasympt}
  \| g- g^S \|_{C^2_{-1-\delta}(\IR^3\setminus B_\sigma(0))}
  := \sup_{\IR^3\setminus B_\sigma(0)} (r^{1+\delta}|g-g^S| +
  r^{2+\delta}|\nabla ^g - \nabla^S| + r^{3+\delta}|\Ric^g -
  \Ric^S |) 
\end{equation}
satisfies $\gasym<\eta$. We let $\Oeta$ denote a constant for which
$\Oeta \leq C\eta$ if $\eta<\eta_0$ is bounded.

The volume element $\rmd V$ of $g$ is a scalar multiple of the
volume element $\rmd V^S$ of $g^S$, that is $\rmd V = h\rmd
V^S$. The asymptotics \eqref{eq:gasympt} imply that $|h|\leq \Oeta
r^{-1-\delta}$. In addition, the scalar curvature $\Scal$ of $g$ satisfies
$|\Scal|\leq \Oeta r^{-3-\delta}$.

Consider a surface $\Sigma\subset \IR^3\setminus B_\sigma(0)$. Let
$\gamma^e, \gamma^S$, and $\gamma$ be the first fundamental forms of
$\Sigma$ induced by $g^e$, $g^S$, and $g$, respectively, $A^e$,
$A^S$, and $A$ the corresponding second fundamental forms, 
$H^e$, $H^S$, and $H$ the mean curvatures and $\Acirc^e$,
$\Acirc^S$, and $\Acirc$ the respective trace free parts of the
second fundamental form.

From the well known transformation behavior for the following
geometric quantities under conformal transformations, and the
asymptotics \eqref{eq:gasympt}, we see:
\begin{lemma}
  \label{lemma:extcurv}
  The normals $\nu^e$, $\nu^S$, and $\nu$ of $\Sigma$ in the metrics
  $g^e$, $g^S$, and $g$ satisfy 
  \begin{gather*}
    \nu^S = \phi^{-2}\nu^e \,,\\
    |\nu^S - \nu| \leq \Oeta r^{-1-\delta} \,,\\
    | \nabla^g\nu^S - \nabla^g\nu| \leq \Oeta r^{-2-\delta}\,.
  \end{gather*}
  The area elements $\rmd\mu^e$, $\rmd\mu^S$, and $\rmd\mu$ satisfy
  \begin{gather*}
    \rmd\mu^S = \phi^4\rmd\mu^e\,,\\
    \rmd\mu - \rmd\mu^S = h\,\rmd\mu \quad \text{with}\quad |h| \leq \Oeta r^{-1-\delta} \,.
  \end{gather*}
  The trace free parts $\Acirc^e,\Acirc^S$, and $\Acirc$ of the second
  fundamental forms satisfy
  \begin{gather*}
    \Acirc^S = \phi^{-2}\Acirc^e\,,  \\
    |\Acirc - \Acirc^S| \leq \Oeta r^{-2-\delta}\,.
  \end{gather*}
  The mean curvatures $H^e,H^S$, and $H$ are related via
  \begin{gather*}
    H^S = \phi^{-2} H^e + 4\phi^{-3}\del_\nu\phi \,,\\
    |H - H^S| \leq \Oeta r^{-2-\delta} \,.
  \end{gather*}
\end{lemma}
To obtain integral estimates for asymptotically decaying quantities,
we cite the following lemma from \cite[Lemma 5.2]{HY:1996}.
\begin{lemma}
  \label{lemma:integral}
  Let $(M,g)$ be $(m,0,\sigma,\eta)$-asymptotically flat, and let $p_0>2$
  be fixed. Then there exists $c(p_0)$, and $r_0=r_0(m,\eta,\sigma)$,
  such that for every hypersurface $\Sigma\subset\IR^3\setminus
  B_{\rmin}(0)$, and every $p>p_0$, the following estimate holds
  \[ \int_\Sigma r^{-p}\dmu \leq c(p_0) \rmin^{2-p}\int_\Sigma
  H^2\dmu\,.\]
  Integration and mean curvature refer to $g$, and $r$ is the
  Euclidean radius.
\end{lemma}
Using lemma \ref{lemma:extcurv} to compare the $L^2$-norms of
$\Acirc$ in the $g$-metric and $\Acirc^S$ in the $g^S$-metric, and
using the conformal invariance of $\|\Acirc^S\|_{L^2(\Sigma,g^S)}$, we obtain
\begin{lemma}
  Let $(M,g)$ be $(m,\delta,\sigma,\eta)$-asymptotically flat. Then there exists
  $r_1=r_1(\eta,\sigma)$, such that for every surface
  $\Sigma\subset\IR^3\setminus B_{\rmin}(0)$ with $\rmin>r_1$, we have
  \begin{multline*}
    \left| \|\Acirc^e\|^2_{L^2(\Sigma,g^e)} - \|\Acirc\|^2_{L^2(\Sigma,g)}
      \right|
      \\ 
      \leq
      \Oeta\rmin^{-1-\delta}\left(\|\Acirc\|^2_{L^2(\Sigma,g)}
      + \|H\|_{L^2(\Sigma)}\|\Acirc\|_{L^2(\Sigma)} +  \rmin^{-1-\delta}\|H\|^2_{L^2(\Sigma)}
      \right) \,.
  \end{multline*}
\end{lemma}
\begin{corollary}
  \label{koro:acirc}
  Let $M$, $g$, $r_1$ and $\Sigma$ be as in the previous lemma. Assume in
  addition that $\|H\|_{L^2(\Sigma)} \leq C'$, then
  \[ \| \Acirc^e \|_{L^2(\Sigma)} \leq C(r_1)\|\Acirc\|_{L^2(\Sigma,g)} +
  C(r_1,C_0,C') \Oeta \rmin^{-1-\delta}\,.
  \]
\end{corollary}
Next we quote a Sobolev-inequality for surfaces contained in
asymptotically flat manifolds. It can be found in \cite[Proposition
  5.4]{HY:1996}. The proof uses the well known Michaels-Simon-Sobolev
inequality in Euclidean space \cite{MS:1973}. 
\begin{proposition}
  \label{prop:sobolev}
  Let $(M,g)$ be $(m,0,\sigma,\eta)$-asymptotically flat. Then there
  is $r_0 = r_0(m,\eta,\sigma)$, and an absolute constant
  $C_\mathrm{sob}$, such that each surface $\Sigma\subset
  M\setminus B_{r_0}(0)$ and each Lipschitz function $f$ on $\Sigma$
  satisfy
  \[ \left( \int_\Sigma |f|^2\,\rmd\mu \right)^{1/2} \leq C_\mathrm{sob}
  \int_\Sigma |\nabla f| + |Hf|\, \rmd\mu\,. \] 
  Using H\"older's inequality, this implies that for all $q\geq 2$
  \begin{equation}
    \label{eq:ms_sob_a}
    \int_\Sigma f^q\dmu \leq 
    C_\mathrm{sob} \left( \int_\Sigma |\nabla u|^\frac{2q}{2+q} +
    |uH|^\frac{2q}{2+q}\right)^\frac{2+q}{2}\,,
  \end{equation}
  and for all $p\geq 1$
  \begin{equation}
    \label{eq:ms_sob_b}
    \left( \int_\Sigma |f|^{2p}\,\rmd\mu \right)^{1/p} \leq
    C_\mathrm{sob}\, p^2\,|\supp f|^{1/p}\int_\Sigma|\nabla f|^2 +
    H^2f^2\,\rmd\mu \,.
  \end{equation}
\end{proposition}

  \section{A priori estimates I}
\label{s:apriori1}
We begin by stating rather general a priori estimates for the geometry
of surfaces. For this, let $\Sigma\subset\IR^3\setminus B_\sigma(0)$
be a surface, and let $g$ be $(m,0,\sigma,\eta)$-asymptotically
flat. Let $\rmin:=\min_\Sigma r$ be the minimum of the Euclidean
radius on $\Sigma$. Assume that on $\Sigma$ the following two
conditions are satisfied:
\begin{gather}
  \label{eq:gradH} 
  \int_\Sigma |\nabla H|^2 \dmu 
  \leq C^K\int_\Sigma r^{-4}|A|^2 + r^{-6}\dmu\,,\\
  \label{eq:A} 
  \int_\Sigma u\,|A|^2 \dmu \leq C^B_0 \int_\Sigma u\,\det A\dmu
  \qquad\text{for all}\qquad 0 \leq u\in C^\infty(\Sigma)\,.
\end{gather}
\begin{remark}
  \begin{mynum}
  \item
    The first condition states that in a certain sense the mean
    curvature is nearly constant. This condition will later be implied
    by the equation by which the mean curvature is prescribed.
  \item
    The second condition means that the surfaces are convex. Indeed,
    on smooth surfaces \eqref{eq:A} implies that $|A|^2\leq \det A$
    pointwise. However, we will need that condition \eqref{eq:A} is
    preserved under $W^{2,p}$-convergence of surfaces.
    Huisken and Yau \cite{HY:1996} are able to replace this
    condition by requiring stability of their CMC surfaces. In the
    present case similar reasoning would work, however, stability is
    not a natural condition for our surfaces.\fillbox
  \end{mynum}
\end{remark}
Condition \eqref{eq:A} implies topological
restrictions, and an estimate on the $L^2$-norm of the mean
curvature. 
\begin{lemma}
  There is $r_0=r_0(m,\eta,\sigma,C_0^B)$, such that every compact
  closed surface $\Sigma$ satisfying \eqref{eq:A} and $\rmin>r_0$
  is diffeomorphic to $S^2$ and satisfies
  \begin{equation}
    \label{eq:A'}
    \int_\Sigma H^2\dmu \leq C(m,\eta,C_0^B)\,. 
  \end{equation}
\end{lemma}
\begin{myproof}
  The Gauss equation \eqref{eq:gauss} implies that the Gauss curvature
  $G$ of $\Sigma$ is given by $G=\det A +\Ric(\nu,\nu) -
  \frac{1}{2}\Scal$. Inserting $u\equiv 1$ into \eqref{eq:A} and
  applying lemmas \ref{lemma:extcurv} and \ref{lemma:integral}, we
  obtain
  \[ \int_\Sigma H^2 \dmu \leq C^B_0 \int_\Sigma G \dmu 
  + C^B_0 \int_\Sigma |\Ric| \dmu \leq C^B_0 \chi(\Sigma) + C^B_0
  \rmin^{-1}\int_\Sigma H^2\dmu\,. \] Here $\chi(\Sigma)$ is the Euler
  characteristic of $\Sigma$. If $\rmin$ is large enough, this gives
  $0\leq \|H\|_{L^2(\Sigma)} \leq C_0^B\chi(\Sigma)$, which implies
  $\chi(\Sigma)\geq 0$. If $\chi(\Sigma)=0$ , i.e. $\Sigma$ is
  a torus, then $\int_\Sigma H^2\dmu=0$ whence
  $\|\Acirc\|_{L^2(\Sigma)}=0$. Using Corollary \ref{koro:acirc} and
  theorem \ref{thm:lm}, we obtain that $\Sigma$ is a sphere, a
  contradiction.
\end{myproof}
\begin{proposition}
  \label{prop:rund}
  Let $(M,g)$ be $(m,0,\sigma,\eta)$-asymptotically flat. Then there
  exists $r_0=r_0(m,\eta,\sigma,C_0^B,C^K)$, such that each closed
  surface $\Sigma$ satisfying \eqref{eq:gradH}, \eqref{eq:A}, and
  $\rmin>r_0$ also satisfies
  \[ \int_\Sigma \left|\nabla|\Acirc|\right|^2 + H^2 |\Acirc|^2\,\rmd\mu 
  \leq C(m,\eta,C_0^B,C^K) \rmin^{-4}\,. \]
\end{proposition}
\begin{myproof}
  We begin by computing
  \[ \Acirc^{\alpha\beta}\Delta\Acirc_{\alpha\beta} = \Acirc^{\alpha\beta}(\Delta A_{\alpha\beta} + \gamma_{\alpha\beta} \Delta H)
  = \Acirc^{\alpha\beta}\Delta A_{\alpha\beta}\,,  \]
  since $\Acirc$ is trace free. By
  \[ 2|\Acirc|\Delta|\Acirc| + 2|\nabla|\Acirc||^2 = \Delta |\Acirc|^2 =
  2\Acirc^{\alpha\beta}\Delta \Acirc_{\alpha\beta} + 2|\nabla \Acirc|^2\,, \]
  and 
  \begin{equation}
    \label{eq:simons_grad}
    |\nabla \Acirc|^2 - |\nabla|\Acirc||^2 \geq 0 \,,
  \end{equation}
  we obtain, using the Simons identity \eqref{eq:simons}
  \begin{equation} 
    \label{eq:simonsa}
    \begin{split}
    |\Acirc|\Delta|\Acirc| 
    &\geq
    \Acirc^{\alpha\beta} \nabla_\alpha \nabla_\beta H 
    + H \Acirc^{\alpha\beta}A_\alpha^\delta A_{\beta\delta} 
    - |A|^2 \Acirc^{\alpha\beta} A_{\alpha\beta} 
    + \Acirc^{\alpha\beta} A_\alpha^\delta R_{\eps\beta\eps\delta}
    \\
    & \phantom{\geq}
    + \Acirc^{\alpha\beta} A^{\delta\eps} R_{\delta\alpha\beta\eps}
    + \Acirc^{\alpha\beta}\nabla_\beta\left(\Ric_{\alpha k}\nu^k\right)
    + \Acirc^{\alpha\beta}
    \nabla^\delta\left(R_{k\alpha\beta\delta}\nu^k\right)\,.
    \end{split}
  \end{equation}
  Integration, and partial integration of $|\Acirc|\Delta|\Acirc|$
  renders
  \begin{equation} 
    \label{eq:simonsb}
    \begin{split}
    \int_\Sigma |\nabla|\Acirc||^2 \dmu
    &\leq \int_\Sigma - \Acirc^{\alpha\beta}\nabla_\alpha \nabla_\beta H
    + |A|^2 \Acirc^{\alpha\beta} A_{\alpha\beta} 
    - H \Acirc^{\alpha\beta}A_\alpha^\delta A_{\beta\delta} \dmu 
    \\
    &\phantom{\leq}
    - \int_\Sigma \Acirc^{\alpha\beta} 
    A_\alpha^\delta R_{\eps\beta\eps\delta}
    + \Acirc^{\alpha\beta}A^{\delta\eps} R_{\delta\alpha\beta\eps}\dmu
    \\
    & \phantom{\leq}
    - \int_\Sigma \Acirc^{\alpha\beta}\left(\nabla_\beta\Ric_{\alpha k}\nu^k\right) +
    \Acirc^{\alpha\beta}\nabla^\delta\left( R_{k\alpha\beta\delta}\nu^k\right) \dmu\,.
    \end{split}
  \end{equation}
  In the first line one computes as follows, and estimates, using
  convexity \eqref{eq:A}
  \begin{equation*}
    \begin{split}
      \int_\Sigma |A|^2 \Acirc^{\alpha\beta} A_{\alpha\beta} - H
      \Acirc^{\alpha\beta}A_\alpha^\gamma A_{\beta\gamma}\dmu
      =
      - 2 \int_\Sigma |\Acirc|^2 \det A \dmu 
      \leq
      - \frac{2}{C^B_0} \int_\Sigma |\Acirc|^2 |A|^2 \dmu \,.
    \end{split}
  \end{equation*}
  To recast the second line of \eqref{eq:simonsb}, recall that in
  three dimensions, the Ricci tensor determines the Riemann tensor:
  \begin{equation}
    \label{eq:riem3d}
    R_{ijkl} = \Ric_{ik} g_{jl} - \Ric_{il}g_{jk} 
    - \Ric_{jk}g_{il} + \Ric_{jl}g_{ik}
    - \frac{1}{2} \Scal \left(g_{ik}g_{jl} - g_{il}g_{jk}\right)\,.
  \end{equation}
  This implies that the second line of \eqref{eq:simonsb} can be expressed as
  \begin{equation*}
    \Acirc^{\alpha\beta} A_\alpha^\delta R_{\eps\beta\eps\delta} 
    + \Acirc^{\alpha\beta} A^{\delta\eps} R_{\delta\alpha\beta\eps}
    = 2\Acirc^{\alpha\beta}\Acirc_\alpha^\delta\Ric_{\beta\delta} - |\Acirc|^2\Ric(\nu,\nu) \,.
  \end{equation*}
  Let $\omega = \Ric(\cdot,\nu)^T$ be the tangential projection of
  $\Ric(\cdot,\nu)$ to $\Sigma$. Then partial integration, the
  Codazzi-equations \eqref{eq:codazzi} and \eqref{eq:riem3d} give for
  the first term of \eqref{eq:simonsb}
  \[ - \int_\Sigma \Acirc^{\alpha\beta} \nabla_\alpha \nabla_\beta
  H\,\rmd\mu 
  = \int_\Sigma \big(\tfrac{1}{2}|\nabla H|^2 + \omega(\nabla H)
  \big)\dmu\,. \]
  In the last line of \eqref{eq:simonsb} we compute, using partial
  integration, \eqref{eq:riem3d}, and \eqref{eq:codazzi}
  \begin{equation*} 
    \int_\Sigma \big(\Acirc^{\alpha\beta}
    \nabla_\beta\big(\Ric_{\alpha k}\nu^k\big) +
    \Acirc^{\alpha\beta}\nabla^\delta\big(R_{k\alpha\beta\delta}\nu^k\big)
    \big)\dmu = - \int_\Sigma 2|\omega|^2 + \omega(\nabla H)\dmu\,.
  \end{equation*}
  Combining these estimates with \eqref{eq:simonsb} and
  $2|\omega(\nabla H)|\leq |\omega|^2 + |\nabla H|^2$, we infer
  \begin{equation*}
    \label{eq:simonsc}
    \int_\Sigma |\nabla|\Acirc||^2 + \tfrac{2}{C^B_0}
    |A|^2|\Acirc|^2\dmu
    \leq
    \int_\Sigma \tfrac{3}{2}|\nabla H|^2 
    + 3 |\omega|^2
    + |\Acirc|^2\Ric(\nu,\nu) -
    \Acirc^{\alpha\beta}\Acirc_\alpha^\delta\Ric_{\beta\delta}\dmu\,. 
  \end{equation*}
  The asymptotics of $g$ imply that $|\Ric|+|\omega|\leq
  C(m,\eta) r^{-3}$. Inserting the estimate
  \eqref{eq:gradH} for $\int_\Sigma|\nabla H|^2\dmu$ into the
  previous estimate, we arrive at
  \begin{equation}
    \label{eq:simonsd}
    \int_\Sigma |\nabla|\Acirc||^2 + \tfrac{2}{C^B_0} |A|^2|\Acirc|^2\,\rmd\mu
    \leq C(m,\eta,C^K)\int_\Sigma \big( r^{-4}|A|^2 + r^{-6}
    + r^{-3}|\Acirc|^2 \big)\dmu\,.
  \end{equation}
  The first term on the right equals $|A|^2 = |\Acirc|^2 +
  \frac{1}{2}H^2$. Using \eqref{eq:A'} we obtain
  \[ \int_\Sigma r^{-4}|A|^2\dmu \leq C(m,\eta,C_0^B)\rmin^{-4} +
  \int_\Sigma r^{-4}|\Acirc|^2\dmu\,.\]
  The third integrand of the right of \eqref{eq:simonsd}
  can be estimated together with the last term of this equation by
  combining the Schwarz inequality with lemma \ref{lemma:integral}    
  \begin{equation*}
    \int_\Sigma r^{-3}|\Acirc|^2 \dmu 
    \leq \int_\Sigma \tfrac{2}{C^B_0} |\Acirc|^{4} + \tfrac{C^B_0}{4r^6}\dmu 
    \leq \frac{2}{C^B_0}\int_\Sigma |\Acirc|^{4} \dmu + C(m,\eta,C_0^B)\rmin^{-4}\,.
  \end{equation*}
  Inserting these estimates into \eqref{eq:simonsd}, and absorbing the
  first term of this equation on the left hand side we obtain the
  assertion of the proposition.
\end{myproof}
\begin{corollary}
  \label{koro:l2_a0}
  Under the additional assumption that
  $(C_1^B)^{-1} R(\Sigma)^{-1} \leq |H|$,
  the previous proposition gives an estimate for the $L^2$-norm
  of $\Acirc$,
  \[ \|\Acirc\|_{L^2(\Sigma)}
  \leq C(m,\eta,C_0^B, C_1^B,C^K)R(\Sigma) \rmin^{-2}\,.\]
\end{corollary}
\begin{corollary}
  \label{koro:grad_a0}
  Under the assumptions of Proposition \ref{prop:rund}, in
  fact
  \begin{equation*}
    \|\nabla \Acirc\|_{L^2(\Sigma)} 
    + \| H \Acirc\|_{L^2(\Sigma)}
    \leq C(m,\eta,C_0^B,C^K)\rmin^{-2}\,.
  \end{equation*}
\end{corollary}
\begin{myproof}
  The proof works by replacing equation \eqref{eq:simons_grad} in the
  proof of proposition \ref{prop:rund} by 
  \begin{equation*}
    |\nabla \Acirc|^2 - |\nabla |\Acirc||^2 \geq \tfrac{1}{17} |\nabla
    \Acirc|^2 - \tfrac{16}{17} (|\omega|^2 + |\nabla H|^2 )\,. 
  \end{equation*}
  This inequality is proved in the same way as a similar inequality
  for $\nabla A$, which is recorded in \cite[Section 2]{SY:1981}. The
  right hand side introduces the desired term $|\nabla \Acirc|^2$, and
  the remaining terms are treated as in the proof of proposition
  \ref{prop:rund}.
\end{myproof}
\begin{corollary}
  Under the assumptions of proposition \ref{prop:rund} and corollary
  \ref{koro:l2_a0}, $\Acirc$, as well as $H$, are controlled in
  the $W^{1,2}$-norm. We therefore have uniform
  estimates for the second fundamental form:
  \begin{align*}
    \|A\|_{L^2(\Sigma)} &\leq
    C(m,\eta,C_0^B,C_1^B,C^K)
    \left(1 + \rmin^{-2}R(\Sigma)\right) \quad\text{and}\\
    \|\nabla A\|_{L^2(\Sigma)} &\leq
    C(m, \eta, C_0^B, C^K)\rmin^{-2}\,.
  \end{align*}
\end{corollary}

  \section{A priori estimates II}
\label{s:apriori2}
This section specializes on surfaces which satisfy the equation 
\begin{equation}
  \label{eq:Heq}
  H\pm P=\const\,.
\end{equation}
We will use theorem \ref{thm:lm} to derive estimates
on the position of such a surface by using the curvature estimates of
the previous section.

As described in the introduction, $P=\tr^\Sigma K = \tr^M K -
K(\nu,\nu)$ is the trace of an extra tensor field $K$ along
$\Sigma$. We will consider data $(M,g,K)$ which are
$(m,\delta,\sigma,\eta)$-asymptotically flat. That is, in addition to
\eqref{eq:gasympt} we have that the weighted norm of $K$ satisfies
\[ \|K\|_{C^1_{-2-\delta}(\IR^3\setminus B_\sigma(0))} 
:=\!\!\sup_{\IR^3\setminus B_\sigma(0)} ( r^{2+\delta}|K| +
r^{3+\delta} |\nabla^g K|) <\eta \,.\] In the sequel, we will consider
either $(m,\delta,\sigma,\eta)$-asymptotically flat data with
$\delta>0$ and arbitrary $\eta<\infty$, or $(m,0,\sigma,\eta)$-
asymptotically flat data with small $\eta\ll 1$.
\begin{remark}
  \label{rem:Heq}
  If $(M,g,K)$ are $(m,\delta,\sigma,\eta)$-asymptotically flat,
  equation \eqref{eq:Heq} implies condition \eqref{eq:gradH}. Indeed
  $|\nabla H|^2 = |\nabla P|^2$ and 
\begin{equation*}
  \nabla^\Sigma P = \nabla^\Sigma  \tr^M K - (\nabla^M_\cdot
 K)(\nu,\nu) - 2K( A(\cdot),\nu)\,,
\end{equation*}
such  that $|\nabla P|^2\leq |\nabla K|^2 + |A|^2|K|^2$. Then
\[ \int_\Sigma |\nabla H|^2 \dmu = \int_\Sigma |\nabla P|^2\dmu \leq
\Kasym^2\int_\Sigma r^{-4-2\delta}|A|^2 +
r^{-6-2\delta}\dmu\,.\fillbox\]
\end{remark}
The results of this section require some additional conditions on the
surfaces:
\begin{gather}
  \label{eq:B1} \tag{A1}
  R(\Sigma) \leq C^A_1 \rmin^q \qquad q<\tfrac{3}{2}\ \text{for}\ 
  \delta>0\quad\text{or}\quad q=1\ \text{for}\ \delta = 0 \,,\\[1mm]
  \label{eq:B2} \tag{A2}
  (C_2^A)^{-1} R(\Sigma)^{-1} \leq H\pm P \,,\\[1mm]
  \label{eq:B3} \tag{A3}
  \int_\Sigma u\,|A|^2 \dmu \leq C^A_3 \int_\Sigma u\,\det A\dmu \qquad\text{for
    all}\qquad 0 \leq u\in C^\infty(\Sigma)\,,\\
  \label{eq:B4} \tag{A4}
  \frac{1}{4\pi R_e(\Sigma)^2}\left| \int_\Sigma \id_\Sigma \dmu^e \right| \leq R_e \,.
\end{gather}
In the sequel, $C^A$ will denote constants which depend only on
$C_1^A, C_2^A$ and $C_3^A$. If $(M,g,K)$ is
$(m,\delta,\sigma,\eta)$-asymptotically flat with $\delta>0$, with
$o(1)$ we denote constants depending on $m$, $C^A$, $\delta$ and
$\eta$, such that $o(1)\to 0$ for $\sigma<\rmin\to\infty$. If
$(M,g,K)$ is $(m,0,\sigma,\eta)$-asymptotically flat, $o(1)$ is such
that, for each $\eps>0$ there is $\eta_0$, and $r_0$ such that
$|o(1)|<\eps$, provided $\eta<\eta_0$ and $\rmin>r_0$. For fixed $m$
and bounded $C^A$, both $r_0$ and $\eta_0$ can be chosen independent of
$C^A$.
\begin{remark}
  \label{rem:conditions}
  \begin{mynum}
  \item 
    Conditions \eqref{eq:B1} and \eqref{eq:B2} allow to compare
    different radius expressions, namely the Euclidean radius $r$, the
    geometric radius $R(\Sigma)$ and the curvature radius given by
    $2/H$. This is necessary, since the curvature estimates of the
    previous section improve with growing $\rmin$, while the estimates
    of DeLellis and M\"uller include the geometric radius
    $R(\Sigma)$. To balance these two radii we use
    \eqref{eq:B1}. Condition \eqref{eq:B2} will be used to apply
    corollary \ref{koro:l2_a0} to obtain $L^2$-estimates on $\Acirc$.
  \item
    Condition \eqref{eq:B4} means that the surface is not far off
    center. We will use this to conclude that the origin is
    contained in the approximating sphere of theorem \ref{thm:lm}.
  \item
    The distinction of the cases $\delta>0$ and $\delta=0$ in
    condition \eqref{eq:B1} is due to the fact that in the proof of
    Proposition \ref{prop:pos}, we can use lemma \ref{lemma:integral}
    only for $\delta>0$.
  \item
    To prove the uniqueness result of Huisken and Yau
    \cite[5.1]{HY:1996} we do not need conditions \eqref{eq:B2} and
    \eqref{eq:B3}. Instead, if we impose stability of the CMC
    surfaces, the estimates $\|\Acirc\|_{L^2(\Sigma)} \leq
    C\rmin^{-1/2}$ and \eqref{eq:A'} can be derived as in
    \cite[5.3]{HY:1996}. Condition \eqref{eq:B1} is slightly stronger
    than what Huisken and Yau need, they only require $q<2$. Using
    stability and \eqref{eq:B1} and \eqref{eq:B4}, only we can prove
    all subsequent estimates.  \fillbox
  \end{mynum}
\end{remark}
The position estimates we will obtain here are formulated in the following:
\begin{proposition}
  \label{prop:pos}
  Let $(M,g,K)$ be $(m,\delta,\sigma,\eta)$-asymptotically flat with $m>0$,
  and let $\Sigma$ be a surface which satisfies \eqref{eq:Heq} and
  \eqref{eq:B1}--\eqref{eq:B4}. Let $R_e$ denote the geometric radius of
  $\Sigma$ and $a$ its center of gravity, both taken with respect to
  the Euclidean metric. Let $S:=S_{R_e}(a)$ denote the Euclidean
  sphere with center $a$ and radius $R_e$. Then there exists a
  conformal parameterization $\psi:S\to(\Sigma,\gamma^e)$, such that
  \begin{gather}
    \label{eq:sup_est}
    {\textstyle \sup_{S}}|\psi-\id_{S}| \leq C(m,C^A) R(\Sigma)^2 \rmin^{-2} \,,\\
    \label{eq:conf_est}
    \|h^2-1\|_{L^2(S)} \leq C(m,C^A) R(\Sigma)^2 \rmin^{-2}\quad\text{and} \\
    \label{eq:norm_est}
    \| N\circ\id_{S} - \nu\circ\psi \|_{L^2(S)}\leq
    C(m,C^A) R(\Sigma)^2 \rmin^{-2} \,.
  \end{gather}
  In addition, the center satisfies the estimate
  \begin{equation} 
    \label{eq:center_est}
    |a|/R_e \leq o(1)\,,
  \end{equation}
  where $o(1)$ is as described at the beginning of this section.
\end{proposition}
\begin{myproof}
  Using \eqref{eq:Heq}, remark \ref{rem:Heq}, and condition
  \eqref{eq:B3}, corollaries \ref{koro:l2_a0} and \ref{koro:acirc} imply
  the following roundness estimates with respect to the Euclidean metric
  \begin{equation}
    \label{eq:rund1}
    \| \Acirc^e\|_{L^2(\Sigma,g^e)} \leq C(m,C^A) R_e \rmin^{-2}
  \,.
  \end{equation}
  Therefore theorem \ref{thm:lm} and the subsequent remarks as
  well as lemma \ref{lemma:extcurv} imply $R_e\leq 2R(\Sigma)$ and
  \eqref{eq:sup_est}--\eqref{eq:norm_est}.
  Condition \eqref{eq:B1} then implies that $R_e\rmin^{-1}\leq
  2C_1^A\rmin^{q-1}$, from which by \eqref{eq:sup_est} 
  \[ |\id_{S}| \geq |\psi| - C(m,C^A)\rmin^{2q-2} \geq r - \frac{1}{2}\rmin
  \geq \frac{1}{2}r \geq \frac{1}{2}\rmin\,, \]
  if $\rmin$ is large enough. Every convex
  combination with $0\leq\lambda\leq 1$ also satisfies
  \begin{equation}
    \label{eq:sph_konvex}
    |\lambda \id_{S_{R_e}(a)} + (1-\lambda)\psi| \geq \frac{1}{2} r\,.
  \end{equation}
  Similar to Huisken and Yau \cite{HY:1996}, we compute for a fixed vector $b\in\IR^3$ with $|b|^e=1$
  \begin{equation}
    \label{eq:pos1}
    0 = (H\pm P) \int_\Sigma g^e(b,\nu^e)\rmd\mu^e 
    = \int_\Sigma H g^e(b,\nu^e)\rmd\mu^e \pm \int_\Sigma P g^e(b,\nu^e)\rmd\mu^e\,.
  \end{equation}
  We estimate using \eqref{eq:B1},
  \begin{equation}
    \label{eq:pos_p_est}
    \left| \int_\Sigma P g^e(b,\nu^e)\rmd\mu^e\right| \leq o(1) \int_\Sigma r^{-2-\delta}
  g^e(b,\nu^e)\rmd\mu^e \leq o(1)\,. 
  \end{equation}
  This follows from lemma \ref{lemma:integral} in the case $\delta>0$,
  and by brute force and \eqref{eq:B1} in the case $\delta=0$. In the
  first term we express $H$ by $H^e$. Using lemma
  \ref{lemma:extcurv} we obtain that the error is  of the order
  $o(1)$, such that  
  \begin{equation}
    \label{eq:posb}
    \left| \int_\Sigma \left(H^e \phi^{-2} +
    4\phi^{-3}\del_{\nu^e}\phi\right)g^e(b,\nu^e)\,\rmd\mu^e \right|
    \leq o(1)\,.
  \end{equation}
  The first variation formula with respect to the Euclidean
  metric gives
  \[ \int_\Sigma H^e\phi^{-2}g^e(b,\nu^e)\rmd\mu^e =
  \int_\Sigma\div_\Sigma^e(\phi^{-2}b)\,\rmd\mu^e =
  -2\int_\Sigma\phi^{-3}g^e(b,\nabla^e\phi)\,\rmd\mu^e\,.\]
  Using $g^e(\nabla^e\phi,b) =g^e(D\phi,b) -
  g^e(b,\nu^e)\del_{\nu^e}\phi$ and $|D\phi|\leq
  C(m)r^{-2}$ gives
  \begin{equation}
    \label{eq:posc}
    \left| \int_\Sigma 6 g^e(b,\nu^e)\del_{\nu^e}\phi\,\rmd\mu^e - \int_\Sigma
    2 g^e(b,D\phi)\,\rmd\mu^e\right| \leq o(1)\,.
  \end{equation}
  Now we will use that $\Sigma$ is approximated by the sphere $S$ as
  described by \eqref{eq:sup_est}--\eqref{eq:norm_est}, and replace the
  integrals of \eqref{eq:posc} by integrals over $S$. For the first term estimate
  \begin{multline*}
    \left| \int_\Sigma g^e(b,D\phi)\,\rmd\mu^e - \int_{S}
    g^e(b,D\phi)\,\rmd\mu^e\right| 
    \\
    \leq 
    \left| \int_{S} (h^2-1)  g^e(b,D\phi\circ\psi)\,\rmd\mu^e\right| 
    + \left| \int_{S}  g^e(b,(D\phi)\circ\psi - D\phi)\,\rmd\mu^e\right|\,.
  \end{multline*}
  Using \eqref{eq:sup_est}, \eqref{eq:conf_est}, \eqref{eq:sph_konvex},
  $|D\phi|\leq C(m)r^{-2}$, $|D^2\phi|\leq C(m) r^{-3}$, and Lemma
  \ref{lemma:integral} we can estimate the error terms
  \begin{multline*}
    \left| \int_{S} (h^2-1)
    g^e(b,D\phi\circ\psi)\,\rmd\mu^e\right|
    \\
    \leq \|h^2-1\|_{L^2(S)}\,\|g^e(b,D\phi\circ\psi)\|_{L^2(S)}
    \leq C \|\Acirc^e\|_{L^2(\Sigma)} R_e \rmin^{-1}\,,
  \end{multline*}
  and
  \begin{multline*}
    \left| \int_{S}
    g^e(b,(D\phi)\circ\psi - D\phi)\,\rmd\mu^e\right| 
    \\
    \leq \sup_S|\psi-\id| \int_{S}\!\!\big(\max_{\lambda\in[0,1]}
    \left|D^2\phi\left(\lambda\id+(1-\lambda)\psi\right)\right|\big)\rmd\mu 
    \leq C \|\Acirc^e\|_{L^2(\Sigma)} R_e \rmin^{-1}\,.
  \end{multline*}
  The second term in \eqref{eq:posc} can be replaced similarly, with
  analogous treatment of the error terms, additionally using
  \eqref{eq:norm_est}. In the end the error is also
  controlled by $C \|\Acirc^e\|_{L^2(\Sigma)} R_e \rmin^{-1}$.
  Therefore, both error terms can be estimated by
  $C(m,C^A)R^2_e\rmin^{-3}$. Using \eqref{eq:B1} gives $R_e^2 \leq
  C^A\rmin^{-q}$, and finally \eqref{eq:posc} implies that
  \begin{equation}
    \label{eq:posd}
    \left| \int_{S} 6 g^e(b,N)\del_N\phi\,\rmd\mu^e - \int_{S}
    2 g^e(b,D\phi)\,\rmd\mu^e\right| 
    \leq o(1) \,. 
  \end{equation}
  Set $b= \frac{a}{|a|}$, and choose coordinates $\varphi$ and
  $\vartheta$ on $S$ such that $g^e(b,N)=\cos\varphi$. Compute
  $D\phi = -\frac{m}{2r^2} \rho$, $N=R_e^{-1}(x-a)$, and $g^e(N,\rho) = R_e
  r^{-1} + r^{-1}|a|\cos\varphi$, where again $\rho=x/r$ is the
  radial direction of $\IR^3$. Inserting this into \eqref{eq:posd} gives
  \begin{equation}
    \label{eq:pose}
    \left| m\int_{S} 3|a|r^{-3} \cos^2\varphi + 2 R_e
    r^{-3}\cos\varphi-|a|r^{-3}\,\rmd\mu^e\right| \leq o(1) \,.
  \end{equation}
  From condition \eqref{eq:B4} we conclude that $|a|\leq R_e$. Using the
  integration formula
  \[ \int_{S_{R_e}(a)} r^{-k}\cos^l\varphi\,\rmd\mu^e = \frac{\pi R_e}{|a|}
  (2R_e|a|)^{-l} \int_{R_e-|a|}^{R_e+|a|} r^{1-k}(r^2-R_e^2-|a|^2)^l\,\rmd r \] 
  we compute the terms in \eqref{eq:pose} and obtain
  \begin{equation}
    \label{eq:posfinal}
    8\pi m |a|/R_e \leq o(1)\,.
  \end{equation}
  Since $m>0$ this implies the last assertion of the proposition.
\end{myproof}
\begin{corollary}
  \label{koro:B1'}
  For each $\eps>0$ we can choose $o(1)$ sufficiently small such that
  \eqref{eq:B1} can be replaced by the stronger assumption
  \begin{equation}
    \label{eq:B1'}
    (1+\eps)^{-1}R(\Sigma)\leq \rmin \leq (1+\eps) R(\Sigma)\,.
  \end{equation}
  In addition \eqref{eq:B4} can be replaced by the assumption
  \begin{equation}
    \label{eq:B4'}
    \frac{1}{4\pi R_e^2} \left|\int_\Sigma \id_\Sigma \dmu^e\right| \leq \eps R_e \,,
  \end{equation}
  provided $\rmin>r_0$, and $r_0=r_0(\eps,m,\sigma,C^A)$ is large enough.
\end{corollary}
\begin{myproof}
  From the position estimates \eqref{eq:sup_est} and \eqref{eq:center_est}
  we obtain for every $p\in S$
  \[ (1- o(1))R_e \leq R_e - |a| \leq |\id_S(p)| \leq |\psi(p)| +
  C(m,C^A)R_e^2\rmin^{-2}\,.\] 
  Since the left hand side is independent of
  $p$, by arranging that $|o(1)|<\eps$ we obtain
  \[ (1-\eps)R_e\leq \rmin + C(m,C^A) R_e^2 \rmin^{-2} \]
  which implies the corollary in view of \eqref{eq:B1}.
\end{myproof}
\begin{corollary}
  \label{koro:B2'}
  Condition \eqref{eq:B2} holds with improved constants. In
  addition the following upper bound is also true
  \begin{equation}
    \label{eq:B2'} 
    (1+\eps)^{-1} R(\Sigma)^{-1} \leq \frac{H}{2} \leq (1+\eps) R(\Sigma)^{-1}\,,
  \end{equation}
  provided $\eta<\eta_0$ is small enough and $\rmin>r_0$ is large enough.
\end{corollary}
\begin{myproof}
  Using the first part of theorem \ref{thm:lm}, the roundness
  estimates \eqref{eq:rund1}, and $|H-H^e| \leq C r^{-2}$ from lemma
  \ref{lemma:extcurv}, we obtain the following estimate for $H$
  \[ \|H - 2/R_e\|^2_{L^2(\Sigma)} \leq C(m,C^A) \rmin^{-2}\,. \]
  By equation \eqref{eq:Heq}, the mean curvature $H$ is nearly
  constant, whence we derive
  \begin{equation*}
    (H\pm P - 2/R_e)^2 |\Sigma| \leq 2\|H-2/R_e\|^2_{L^2(\Sigma)} + 2 \|P\|^2_{L^2(\Sigma)} \leq
    C(m,C^A) \rmin^{-2}\,. 
  \end{equation*}
  In view of $|P|\leq C(\|K\|_{C^1_{-2}})r^{-2}$, this implies
  $|H-2/R_e| \leq C(m,C^A)\rmin^{-2}$, which gives the assertion of
  the corollary.
\end{myproof}
We now take a closer look at those terms in the proof of proposition
\ref{prop:rund} which came from the geometry of $M$. The Ricci tensor
of the Schwarzschild metric, when restricted to a centered coordinate
sphere for example splits orthogonally into a positive tangential part
$(\Ric^S)^T = mr^{-3}\phi^{-6}\gamma^S\geq 0$ and a negative normal
part $\Ric^S(\nu,\nu) = -2mr^{-3}\phi^{-6} \leq 0$, the mixed term
$\omega^S$ vanishes. We now combine the estimates of proposition
\ref{prop:pos} to estimate the analogous terms on $\Sigma$.
\begin{proposition}
  \label{prop:center_curv}
  Let $\Sigma$ be as in proposition \ref{prop:pos}, then for
  $\rmin>r_0(m,\sigma,C^A)$ large enough, we have
  \begin{align*}
    \|\nu - \phi^{-2}\rho\|^2_{L^2(\Sigma,g)} 
    &\leq 
    o(1)\rmin^2 + C(m,C^A) \,, 
    \\
    \|\Ric(\nu,\nu) - \phi^{-4}\Ric^S(\rho,\rho)\|^2_{L^2(\Sigma,g)} 
    &\leq
    o(1)\rmin^{-4} + C(m,C^A)\rmin^{-6}\,, 
    \\
    \|\omega\|_{L^2(\Sigma,g)}^2 
    &\leq 
    o(1)\rmin^{-4} + C(m,C^A)\rmin^{-6} \,,  
    \\
    \|\Ric^T - P^S_{\phi^{-2}\rho}\Ric^S\|_{L^2(\Sigma,g)}^2 
    &\leq 
    o(1)\rmin^{-4} + C(m,C^A)\rmin^{-6} \,,
  \end{align*}
  where $\rho = x/r$ is the radial direction of $\IR^3$,
  $P^S_{\phi^{-2}\rho}\Ric^S$ is the $g^S$-orthogonal projection of
  the Ricci tensors of $g^S$ onto the subspace of the tangential space
  of $M$ which is $g^S$-orthogonal to ${\phi^{-2}\rho}$, and $o(1)$ is
  as described at the beginning of the section.
\end{proposition}
\begin{myproof}
  From lemma \ref{lemma:extcurv} and corollary \ref{koro:B1'} we
  derive
  \begin{equation}
    \label{eq:nuabsch0}
    \|\nu-\phi^{-2}\rho\|_{L^2(\Sigma,g)}^2 \leq c \|\nu^e-\rho\|^2_{L^2(\Sigma,g^e)} + o(1)\,.
  \end{equation}
  Now we use proposition \ref{prop:pos} to obtain a sphere $S=S_{R_e}(a)$ and a
  conformal parameterization $\psi:S\to \Sigma$ satisfying the
  estimates \eqref{eq:sup_est}--\eqref{eq:center_est}. From the
  estimate on the center $a$, we compute for the difference of the
  Euclidean normal $N = (x-a)/R_e$ and the radial direction $\rho=x/r$ that
  \begin{equation*}
    \label{eq:nuabsch1}
    |N-\rho|_{g^e} \leq (|R^e-r|_{g^e}+|a|)/R^e \leq 2|a|/R^e \leq o(1) \,.
  \end{equation*}
  Using \eqref{eq:sup_est}--\eqref{eq:norm_est} we estimate 
  \begin{equation*}
    \label{eq:nuabsch2}
    \begin{split}
      | \rho \circ \psi(x) - \rho(x) |_{g^e} 
      &\leq 
      \big(\sup_{\lambda\in [0,1]}
      |D\rho(\lambda x - (1-\lambda)\psi(x)|_{g^e}\big)
      |\psi(x)-x|_{g^e}
      \\
      &\leq
      C^A \|\Acirc^e\|_{L^2(\Sigma,g^e)}\,,
    \end{split}
  \end{equation*}
  and
  \begin{equation*}
    \label{eq:nuabsch3}
    \int_\Sigma |\nu^e-\rho|_{g^e}^2\,\rmd\mu^e 
    = \int_S h^{-2} |\nu^e\circ\psi -\rho\circ\psi|_{g^e}^2\rmd\mu^e 
    \leq C \int_S|\nu^e\circ\psi -\rho\circ\psi|_{g^e}^2 \rmd\mu^e\,.
  \end{equation*}     
  By the triangle inequality and the previous inequalities we obtain
  \begin{equation*}
    \label{eq:nuabsch4}
    \begin{split}
      \|\nu^e\circ\psi -\rho\circ\psi\|_{L^2(S,g^e)}
      &\leq \|\nu^e\circ\psi-N\|_{L^2(S)} +
      \|N-\rho\|_{L^2(S)} + \|\rho-\rho\circ\psi\|_{L^2(S)}
      \\
      &\leq o(1)\rmin + C(m,C^A)\,.
    \end{split}
  \end{equation*}
  This implies the first inequality of the proposition in view of \eqref{eq:nuabsch0}.
  The second inequality now easily follows, since 
  \begin{equation*}
    \begin{split}
      &\|\Ric(\nu,\nu) -
      \phi^{-4}\Ric^S(\rho,\rho)\|_{L^2(\Sigma)}^2
      \\      
      &\quad\leq 
      \|\Ric^S-\Ric\|^2_{L^2(\Sigma)} +
      \sup_{\Sigma}|\Ric^S|^2\|\nu-\phi^{-2}\rho\|^2_{L^2(\Sigma)}
      \\
      &\quad\leq 
      C(m, C^A)\rmin^{-6}\big(1 +
      \|\nu-\phi^{-2}\rho\|^2_{L^2(\Sigma)} \big)\,.
    \end{split}
  \end{equation*}
  For the third inequality, observe that by a similar computation
  \begin{equation*}
    \|\Ric(\nu,\cdot) - \Ric^S(\phi^{-2}\rho,\cdot)\|_{L^2(\Sigma)}^2 \leq C(m,
    C^A)\rmin^{-6}\big(1 + \|\nu-\phi^{-2}\rho\|^2_{L^2(\Sigma)} \big)\,,
  \end{equation*}
  such that only the difference of the projections of
  $\Ric^S(\cdot,\phi^{-2}\rho)$ to the subspaces $g$-or\-thog\-o\-nal to
  $\nu$ and $g^S$-orthogonal to $\phi^{-2}\rho$ have to be
  estimated. Note that the latter projection is zero. To estimate the
  difference, write
  \begin{equation*}
    P^g_\nu \Ric^S(\cdot,\phi^{-2}\rho) =  \Ric^S(\cdot,\phi^{-2}\rho) -
    g(\cdot,\nu)\Ric^S(\nu,\phi^{-2}\rho)\,,\\
  \end{equation*}
  where $P_\nu^g$ is the $g$-orthogonal projection on the $g$
  orthogonal complement of $\nu$, and
  \begin{equation*}
    P^S_{\phi^{-2}\rho}\Ric^S(\cdot,\phi^{-2}\rho) = \Ric^S(\cdot,\phi^{-2}\rho)
    - g^S(\cdot,\phi^{-2}\rho)\Ric^S(\phi^{-2}\rho,\phi^{-2}\rho)\,.
  \end{equation*}
  Therefore the third estimate of the proposition follows as before.
  The last estimate can be obtained using a similar computation.
\end{myproof}
We can now improve the roundness estimates of proposition
\ref{prop:rund}.
\begin{proposition}
  \label{prop:rund2}
  Let $(M,g,K)$ be $(m,\delta,\sigma,\eta)$-asymptotically flat. Then there
  exist a constant $C(m,C^A)$ and $r_0=r_0(m,\sigma,C^A)$, such that for
  all surfaces $\Sigma$ satisfying \eqref{eq:Heq}, conditions
  \eqref{eq:B1}--\eqref{eq:B4}, and $\rmin>r_0$, the following estimate
  holds
  \[ \int_\Sigma | \nabla|\Acirc||^2 + H^2|\Acirc|^2 \,\rmd\mu \leq
  o(1) \rmin^{-4} + C(m,C^A)\rmin^{-6}\,. \]
\end{proposition}
\begin{myproof}
  We use the Simons identity as in the proof of proposition \ref{prop:rund}
  \begin{equation*}
    \int_\Sigma |\nabla|\Acirc||^2 + \tfrac{2}{C^A_0}
    H^2|\Acirc|^2\dmu  
    \leq 
    \int_\Sigma (\tfrac{3}{2}|\nabla H|^2 
    + 3|\omega|^2
    + |\Acirc|^2\Ric(\nu,\nu) 
    - \Acirc^{\alpha\beta}\Acirc_\alpha^\delta\Ric_{\beta\delta}
    \dmu\,.
  \end{equation*}
  By Remark \ref{rem:Heq} we have $|\nabla H|^2 \leq o(1)(r^{-4}|A|^2
  + r^{-6})$. We further proceed as in the proof of proposition \ref{prop:rund}
  but now estimate the resulting terms using proposition
  \ref{prop:center_curv}. For example with $\Ric^S(\rho,\rho)\leq 0$
  and the Schwarz inequality we derive
  \begin{equation*}
    \begin{split}
      \int_\Sigma |\Acirc|^2 \Ric(\nu,\nu) \,\rmd\mu
      &\leq
      \|\Acirc\|_{L^4(\Sigma)}^{2}\|\Ric(\nu,\nu) -
      \phi^{-4}\Ric^S(\rho,\rho)\|_{L^2(\Sigma)}
      \\ &\leq
      o(1)\rmin^{-5} + C(m,C^A)\rmin^{-6}\,.
    \end{split}
  \end{equation*}
  Here we used the Sobolev inequality from proposition
  \ref{prop:sobolev} together with proposition \ref{prop:rund} and
   corollary \ref{koro:B1'}, to estimate the $L^4$-norm of $\Acirc$
  \[ \|\Acirc\|_{L^4(\Sigma)}^4 \leq C(m,C^A)|\Sigma|\rmin^{-8}\leq
  C(m,C^A)\rmin^{-6}\,. \] 
  The estimates for the other terms are obvious.  
\end{myproof}
Our next step is to prove $\sup$-estimates for $\Acirc$ using a Stampaccia
iteration.
\begin{proposition}
  \label{prop:rund_sup}
  Let $\Sigma$ be as in proposition \ref{prop:rund}, then for each
  $\eps>0$ there exists $r_0=r_0(m,\sigma,C^A)$ and a constant
  $C(\eps,m,C^A)$, such that if $\rmin\geq r_0$
  \[ \sup_\Sigma | \Acirc| 
  \leq C(\eps,m,C^A)\big(o(1)\rmin^{-2} + \rmin^{-3+\eps}\big)\,. \]
\end{proposition}
\begin{myproof}
  Let $u:=|\Acirc|$, and $u_k:=\max(u-k,0)$ for all $k\geq 0$. Let
  $A(k):= \{x\in\Sigma: u_k>0 \}$. Let $p>1$, and multiply equation
  \eqref{eq:simonsa} with $u_k^p$ and integrate. Partial integration,
  proceeding as in proposition \ref{prop:rund}, and using the Schwarz
  inequality to absorb all gradient terms on the left hand side gives 
  \begin{multline}
    \label{eq:stamp1}
    \IAK p u_k^{p-1}u|\nabla u|^2 + u_k^p |\nabla u|^2 + C^A_3 u_k^p
    u^2 H^2\rmd\mu 
    \\ 
    \leq c(p) \IAK u_k^{p-1} u |\nabla H|^2  + u_k^p|\nabla H|^2 + u_k^p|\omega|^2 
    + u_k^p u^2 |\Ric| + u_k^{p-1} u |\omega|^2\rmd\mu\,.
  \end{multline}
  We have the bounds $|\Ric|+|\omega|\leq C(m) r^{-3}$, and remark \ref{rem:Heq} and corollary
  \ref{koro:B2'} imply that
  $|\nabla H|^2 \leq  o(1) (r^{-6} + r^{-4} u^2)$. Equation
  \eqref{eq:stamp1} therefore gives
  \begin{equation}
    \begin{split}
      \label{eq:stamp2}
      &\IAK p u_k^{p-1}u|\nabla u|^2 + u_k^p |\nabla u|^2 + C^A_3 u_k^p
      u^2 H^2\,\rmd\mu
      \\ 
      &\quad\leq
      C(m,C^A) \IAK u_k^p r^{-6} + u_k^{p-1} u r^{-6} + u_k^p
      u^2 r^{-3} \dmu\,. 
    \end{split}
  \end{equation}
  Using the Sobolev inequality \eqref{eq:ms_sob_b}, proposition
  \ref{prop:rund}, $u_k\leq u$, and $\nabla u_k=\nabla u$ on $A(k)$,
  we infer that for all $1<q<\infty$
  \begin{equation*} 
    \begin{split}
      \IAK u_k^p \rmd\mu 
      &\leq C(q,m,C^A)|A(k)| (o(1)\rmin^{-2p} + \rmin^{-3p}) 
      \\
      &\leq C(q,m,C^A)|A(k)|^{1-1/q} |\Sigma|^{1/q} (o(1)\rmin^{-2p} +
      \rmin^{-3p})\,.
    \end{split}
  \end{equation*}
  We proceed estimating the second term on the right hand side of
  \eqref{eq:stamp2}. We use the Sobolev inequality \eqref{eq:ms_sob_a}
  to conclude that
  \begin{multline*}
      \IAK u_k^{p-1} u\dmu 
      = \IAK (u_k u^\frac{1}{p-1})^{p-1} 
      \\
      \leq C(p)\left(\IAK\, \big|\nabla
      (u_k u^\frac{1}{p-1})\big|^\frac{2(p-1)}{p+1} + |(u_k u^\frac{1}{p-1})
      H|^\frac{2(p-1)}{p+1} \dmu\right)^\frac{p+1}{2}
  \end{multline*}
  Since $|\nabla (u_k u^\frac{1}{p-1})| \leq c(p) u^\frac{1}{p-1}|\nabla
  u|$, we estimate the first term, using H\"older,
  \begin{equation*}
    \begin{split}
      \IAK \, \big|\nabla (u_k u^\frac{1}{p-1})\big|^\frac{2(p-1)}{p+1} \dmu
      &\leq 
      c(p) \IAK |\nabla u|^\frac{2(p-1)}{p+1} u^\frac{2}{p+1}\dmu 
      \\
      &\leq 
      c(p) \left(\IAK |\nabla u|^2\dmu\right)^\frac{p-1}{p+1}\left(\IAK
      u\dmu\right)^\frac{2}{p+1}\,.
    \end{split}
  \end{equation*}
  Similary, for the second term
  \begin{equation*}
    \IAK |(u_k u^\frac{1}{p-1}) H|^\frac{2(p-1)}{p+1} \dmu 
    \leq
    \left( \IAK u^2 H^2 \dmu\right)^\frac{p-1}{p+1}\left(\IAK
    u\dmu\right)^\frac{2}{p+1}\,.   
  \end{equation*}
  Combining these we get that
  \begin{equation*}
    \IAK u_k^{p-1} u\dmu \leq C(p)\left(\IAK u\dmu\right) \left(\IAK
    |\nabla u|^2 + H^2 u^2 \dmu \right)^\frac{p-1}{2}\,.
  \end{equation*}
  Observe that for any $0<q<\infty$, by an application of the
  H\"older inequality and the Soboloev inequality \eqref{eq:ms_sob_b} 
  \begin{equation*}
    \IAK u \dmu 
    \leq 
    \left(\int_\Sigma u^q \dmu \right)^\frac{1}{q}
    |A(k)|^\frac{q-1}{q}
    \leq C(q) |A(k)|^\frac{q-1}{q} |\Sigma|^\frac{1}{q}
    \left(\int_\Sigma  |\nabla u|^2 + H^2 u^2 \dmu\right)^\frac{1}{2}
  \end{equation*}
  In view of proposition \ref{prop:rund2}, this yields that for all
  $1<q<\infty$
  \begin{equation*}
    \IAK u_k^{p-1} u \dmu \leq C(p,q,m,C^A) |A(k)|^\frac{q-1}{q}
    |\Sigma|^\frac{1}{q} (o(1)\rmin^{-2p} + \rmin^{-3p})\,.
  \end{equation*}
  A similar treatment of the last term in equation \eqref{eq:stamp2}
  gives that
  \begin{equation*}
    \IAK u_k^p u^2 \leq C(p,q,m,C^A) |A(k)|^\frac{q-1}{q}
    |\Sigma|^\frac{1}{q} (o(1)\rmin^{-2p-2} + \rmin^{-3p-4})\,.
  \end{equation*}
  Thus we infer that \eqref{eq:stamp2} yields 
  \begin{multline*}
    \IAK p u_k^{p-1}u|\nabla u|^2 + u_k^p |\nabla u|^2 
    + C^A_3 u_k^p u^2 H^2\dmu 
    \\
    \leq C(p,q,m,C^A) |A(k)|^\frac{q-1}{q}|\Sigma|^\frac{1}{q} (o(1)\rmin^{-2p-6} + \rmin^{-3p-6})\,.
  \end{multline*}
  Let $f:=u_k^{p/2+1}$, then the above etimate is equivalent to
  \begin{equation*}
    \IAK |\nabla f|^2 + H^2 f^2 \dmu
    \leq 
    C(p,q,m,C^A) |A(k)|^\frac{q-1}{q}|\Sigma|^\frac{1}{q}(o(1)\rmin^{-2p-6} + \rmin^{-3p-6})\,.
  \end{equation*}  
  Using the Sobolev inequality \eqref{eq:ms_sob_b} to estimate
  $\int_\Sigma f^2\dmu$, and reexpressing this in terms of
  $f^2=u_k^{p+2}$, we obtain the iteration inequality
  \begin{equation*}
    \begin{split}
      |h-k|^{p+2}|A(h)| 
      &\leq 
      \int_{A(h)}\!\!\! u_k^{p+2} \dmu 
      \leq 
      \IAK u_k^{p+2} \dmu
      \\
      &\leq 
      C(p,q,m,C^A)|A(k)|^{2-\frac{1}{q}}|\Sigma|^\frac{1}{q}(o(1)\rmin^{-2p-6}
      + \rmin^{-3p-6})\,.
    \end{split}
  \end{equation*}
  By \cite[Lemma 4.1]{Stampaccia:1966}, this iteration inequality
  implies that $|A(d)|=0$ for $d\geq d_0$ with 
  \begin{equation*}
    d_0^{p+2} \leq  C(p,q,m,C^A)(o(1)\rmin^{-2p-6} +
    \rmin^{-3p-6})|\Sigma|^\frac{1}{q}|A(0)|^{1-\frac{1}{q}}
  \end{equation*}
  As $|A(0)| \leq |\Sigma$ we see that we can fix any
  $1<q<\infty$. Corollary \ref{koro:B1'} implies that
  $|\Sigma|\leq C(m,C^A)\rmin^2$. Therefore we can estimate
  \begin{equation*}
    d_0 \leq  C(p,m,C^A) \big(o(1)\rmin^{-2} + \rmin^{-3+2/(p+2)}\big)\,.
  \end{equation*}
  This yields
  \begin{equation*}
    \sup_\Sigma |\Acirc| \leq C(p,m,C^A) \big(o(1)\rmin^{-2} +
    \rmin^{-3+2/(p+2)}\big)\,.
  \end{equation*}
  Thus the claimed estimate follows provided $p$ is large enough.
\end{myproof}
We now have a $\sup$-estimate for $A=\nabla\nu$. This can be combined
with the $L^2$-estimates for $|\nu-\phi^{-2}\rho|$ to prove a
$\sup$-estimate for this expression.
\begin{proposition}
  \label{prop:norm_sup}
  Let $\Sigma$ be as in proposition \ref{prop:pos} such that in
  particular $\Sigma$ satisfies 
  $\|\nu-\phi^{-2}\rho\|_{L^2(\Sigma)} \leq o(1)\rmin + C(m,C^A)$ and
  $|A| \leq C(m,C^A)\rmin^{-1}$. Then there exists
  $r_0=r_0(m,\sigma,C^A)$ such that 
  \[ \sup |\nu-\phi^{-2}\rho| \leq o(1)+ C(m,C^A)\rmin^{-2/3} \]
  provided $o(1)$ is small enough, and $\rmin>r_0$.
\end{proposition}
 \begin{myproof}
   From the above assumptions, $|\nabla (\nu-\phi^{-2}\rho)| \leq
   C(m,C^A)\rmin^{-1}$. Therefore $f:=|\nu-\phi^{-2}\rho|^2$ satisfies
   \[ |\nabla f| =
   \left|g(\nabla(\nu-\phi^{-2}\rho),\nu-\phi^{-2}\rho)\right|
   \leq C(m,C^A)\rmin^{-1}\,, \]
   provided $r_0$ is large enough. Assume there exists $p_0\in\Sigma$
   such that for $M>0$ the inequality $f(p_0) \geq 2M
   (o(1)+\rmin^{-1})^{2/3}$ holds. Let $B:=\{p\in\Sigma:|p-p_0|\leq M
   (o(1)+\rmin^{-1})^{2/3}C(m,C^A)^{-1}\rmin \}$. Then for all $p\in
   B$ we have that $f(p)\geq M (o(1)+\rmin^{-1})^{2/3}$, which implies
   that
   \[ \int_\Sigma f \,\rmd\mu
   \geq \int_B f \rmd\mu \geq C \frac{M^3}{C(m,C^A)^2}(o(1)\rmin +1)^2\,,
   \]
   where we used that $|B|\geq C M^2
   (\eps+\rmin^{-1})^{4/3}C(m,C^A)^{-2}\rmin^2$. This follows from
   the estimate on the conformal factor of $\psi: S\to \Sigma$ from
   theorem \ref{thm:lm}, if $\eps$ and $\rmin^{-1}$ are small
   enough. If $M$ is large enough, this is a contradiction.
 \end{myproof}
\begin{corollary}
  \label{koro:gradient}
  In the same way we obtain an estimate $\sup_\Sigma|\nu^e-\rho|\leq
  o(1) + C(m,C^A)\rmin^{-2/3}$, and therefore $\int_\Sigma
  g^e(\nu^e,\rho)\geq \frac{1}{2}$, if $o(1)$ is small
  enough. Hence $\Sigma$ is globally a graph over $S^2$,
  i.e. there is a function $u\in C^\infty(S^2)$ such that
  \[ \Sigma = \{ u(p)p : p\in S^2\subset\IR^3 \}\,.\]
\end{corollary}
\begin{corollary}
  \label{koro:center_curv_sup}
  Surfaces $\Sigma$ as in proposition \ref{prop:pos} satisfy
  \begin{equation*}
    |\Ric(\nu,\nu) + 2mr^{-3}| 
    \leq o(1)\rmin^{-3} + C(m,C^A)\rmin^{-3-2/3}\,.
  \end{equation*}
\end{corollary}
This enables us to precisely compute the curvature of $\Sigma$
taken with respect to $g$.
\begin{theorem}
  \label{thm:apriori}
  Let $\Sigma$ be as in proposition \ref{prop:rund}. Let $R_e=
  \sqrt{|\Sigma|^e/4\pi}$ be its Euclidean geometric radius, and define
  $\bar\phi=1+ \frac{m}{2R_e}$ and $\bar H=\frac{2}{\bar\phi^2 R_e} -
  \frac{2m}{\bar\phi^3 R_e^2}$. Then there exist $r_0=r_0(m,\sigma,C^A)$
  and $C(m,C^A)$, such that if $\rmin>r_0$ the following
  estimates hold:
  \begin{align*}
    \sup_\Sigma |H-\bar H| 
    &\leq
    o(1)\rmin^{-2} +
    C(m,C^A)\rmin^{-2-2/3}\,, 
    \\
    \sup_\Sigma\left| \det A - \bar H^2/4\right| 
    &\leq
    o(1)\rmin^{-3} + C(m,C^A)\rmin^{-3-2/3}\,,
    \\
    \sup_\Sigma\left| G - \bar H^2/4 - 2m/R_e^3 \right|
    &\leq
    o(1)\rmin^{-3} + C(m,C^A)\rmin^{-3-2/3}\,.
  \end{align*}
  Here $G = \det A - \Ric(\nu,\nu) + \frac{1}{2}\Scal$ is the
  Gauss-curvature of $\Sigma$.
\end{theorem}
\begin{myproof}
  From Proposition \ref{prop:pos} we obtain an approximating
  sphere $S=S_{R_e}(a)$ and a conformal map $\psi:S\to\Sigma$ which
  satisfies \eqref{eq:sup_est}--\eqref{eq:center_est}. We compare
  $\Sigma$ with the \emph{centered} sphere $\bar S=S_{R_e}(0)$ and
  consider the map $\xi:\bar S \to \Sigma: x\mapsto \psi(x+a)$. From
  \eqref{eq:sup_est} and \eqref{eq:center_est} we obtain that
  \[\sup_\Sigma|r-R^e| 
  = \sup_{\bar S} \big| \left|\xi(x)\right| - \left|x\right| \big|
  \leq o(1)\rmin + C(m,C^A)\,, \] 
  which in particular implies that $|\rmin - R^e| \leq o(1)\rmin +
  C(m,C^A)$. In addition $\sup_\Sigma |\bar\phi -\phi| \leq
  o(1)\rmin^{-1} + C(m,C^A)\rmin^{-2}$ as well as $\sup_\Sigma
  |\bar\phi^{-2} -\phi^{-2}| + \sup_\Sigma |\bar\phi^{-3} -\phi^{-3}|
  \leq o(1)\rmin^{-1} + C(m,C^A)\rmin^{-2}$. Take a point $x\in\bar S$,
  and let $\nu^e$ be the Euclidean normal to $\Sigma$. Estimate
  \begin{equation*}
    \begin{split}
      &| D_\rho(x)\phi(x) - D_{\nu^e(\xi(x))}\phi(\xi(x)) |
      \\
      &\quad 
      \leq 
      | D_{\rho(x)}\phi(x) - D_{\rho(x)}\phi(\xi(x))|
      + | D_{\rho(x)}\phi(\xi(x)) - D_{\rho(\xi(x))}\phi(\xi(x))| 
      \\
      &\quad\phantom{\leq} 
      + | D_{\rho(\xi(x))}\phi(\xi(x)) -
      D_{\nu^e(\xi(x))}\phi(\xi(x))| 
      \\
      &\quad
      \leq 
      o(1)\rmin^{-2} + C(m,C^A)\rmin^{-2-2/3}\,.
    \end{split}
  \end{equation*}
  The $L^2$-norm of $H-\bar H$ can then be estimated by using lemma
  \ref{lemma:extcurv} to replace $H$ by $H^e$, and estimating $\|H^e -
  2/R_e^2\|_{L^2(\Sigma)}$ by taking the trace of \eqref{eq:lm_Aest}.
  \begin{equation*}
    \begin{split}
      \int_\Sigma\!|H-\bar H|^2\,\rmd\mu^e
      &\leq 
      \int_\Sigma \left| \frac{H^e}{\phi^2} + 4\frac{D_{\nu^e}\phi}{\phi^3} -
      \frac{2}{\bar\phi^2 R_e} +
      \frac{2m}{\bar\phi^3 R_e^2}\right|^2\,\rmd\mu^e + o(1)\rmin^{-3}
      \\
      &\leq 
      o(1)\rmin^{-2} + C(m,C^A)\rmin^{-2-4/3}\,.
    \end{split}
  \end{equation*}
  Proceeding as in the proof of corollary \ref{koro:B2'} we obtain the
  asserted $\sup$-estimate on $H-\bar H$ by using \eqref{eq:Heq} and
  $|P|\leq o(1)\rmin^{-2}$.
  
  The $\sup$-estimates on $\Acirc$ of proposition \ref{prop:rund_sup}
  imply that $A$ on $\Sigma$ satisfies
  \begin{equation} 
    \label{eq:A_est_sup}
    \left| A - \frac{1}{2}\bar H\Id \right| \leq \left| A -\frac{1}{2} H\Id \right| + |H-\bar H| 
    \leq o(1)\rmin^{-2} + C(m,C^A)\rmin^{-2-2/3}\,,
  \end{equation}
  which implies the second assertion of the theorem. Corollary
  \ref{koro:center_curv_sup} gives that
  \[ \left|\Ric(\nu,\nu)+2m/R_e^3\right| \leq o(1)\rmin^{-3} +
  C(m,C^A)\rmin^{-3-2/3}\,,\] which, in view of $|\Scal|\leq
  o(1)\rmin^{-3}$, equation \eqref{eq:A_est_sup}, and the Gauss equation
  $G=\det A - \Ric(\nu,\nu) + \frac{1}{2}\Scal$, implies the last
  assertion.
\end{myproof}

  \section{The linearization of the operator \CHP}
\label{s:linearized}
In this section we will examine the linearization of the operator $\CHP$
which assigns the function $H\pm P$ to a surface. We will prove that
this linearization is invertible, whence we can apply the implicit function
theorem in section~\ref{s:foliation} to find surfaces with $H\pm P
=\const$. We begin by computing the linearization. For this let
$\Sigma\subset M$ be a closed surface. In a neighborhood of $\Sigma$ we
introduce Gaussian normal coordinates $y:\Sigma\times(-\eps,\eps)\to
M$, such that $y(\cdot,0)=\id_\Sigma$, and $\del y/\del t =
\nu_{\Sigma_t}$, with $\Sigma_t = y(\Sigma,t)$. For a function $f\in
C^\infty(\Sigma)$ with $|f|\leq\eps$ define the graph of $f$ over
$\Sigma$ as
\[ \graph(f) := \{ y(p,f(p)): p\in\Sigma \}\,. \]
Let $\CH:C^\infty(\Sigma)\to C^\infty(\Sigma)$ be the operator, which
assigns to a function $f$ the mean curvature $\CH(f)$ of
$\graph(f)$, and let $\CP:C^\infty(\Sigma)\to C^\infty(\Sigma)$
be the operator which assigns to a function $f$ the function
$P=\tr^{\graph(f)}K$ evaluated on $\graph(f)$. To compute the
linearization of $\CH\pm\CP$ at $f=0$, we need the following lemma:
\begin{lemma}
  \label{lemma:evol_eq}
  Let $\Sigma\subset M$ be a surface, and $F:\Sigma\times
  (-\eps,\eps)\rightarrow M$ a variation of $\Sigma$,
  with $F(\cdot,0)=\id_\Sigma$. If $F$ is normal to $\Sigma$, i.e. 
  $\ddeval{F}{t}{t=0} = f\nu$ for $f\in C^\infty(\Sigma)$, then
  \begin{gather*}
    \DDeval{H}{t}{t=0} =  -\Delta^\Sigma f - f \left( |A|^2 + \Ric(\nu,\nu) \right) \,,\\
    \DDeval{P}{t}{t=0} = f \left(\nabla^M_\nu \tr^M K - \nabla^M_\nu K(\nu,\nu)\right) +
    2 K(\nabla^\Sigma f,\nu)\,.
  \end{gather*}
  Here $A$ is the second fundamental form, $H$ the mean curvature,
  and $\nu$ the normal of $\Sigma$. The covariant derivative of $M$
  is denoted by $\nabla^M$ and that of $\Sigma$ by $\nabla^\Sigma$.
\end{lemma}
\begin{myproof}
  The first equation is well known. It can be found in \cite[Appendix
  A]{Bray:1997}. The second immediately follows from $P=\tr^M
  K-K(\nu,\nu)$ and $\DDeval{\nu}{t}{t=0}= - \nabla^\Sigma f$.
\end{myproof}
Lemma \ref{lemma:evol_eq} implies that the linearization $\LHP$ of $\CHP$ is given by 
\begin{equation}
  \label{eq:lin1}
  \LHP\!\!f = -\Delta f - f\left(|A|^2+\Ric(\nu,\nu) \pm \nabla^M_\nu\!
    K(\nu,\nu) \mp \nabla^M_\nu\!\!\tr K\right) \pm
  2 K(\nabla^\Sigma\! f,\nu) \,.
\end{equation}
To obtain a form which is easier to handle, we multiply this by $f$ and integrate
by parts. 
\begin{proposition}
  \label{prop:lin_pi}
  Let $f\in C^\infty(\Sigma)$, then
  \begin{equation*}
    \label{eq:lin2}
    \begin{split}
      \int_\Sigma f \LHP f \dmu 
      &= \int_\Sigma |\nabla f|^2 - f^2 \big(8\pi
      \left(\mu \mp J(\nu)\right) + \tfrac{1}{2} \big| (K^T)^\circ \pm
      \Acirc\big|^2 + |\theta|^2\big)
      \\ 
      &\quad -\tfrac{1}{2}f^2\big(\tfrac{1}{2} (H\pm P)^2 + \left(H\mp K(\nu,\nu)\right)^2 - 
      (\tr K)^2 - 2G \big)\rmd\mu\,.\notag
    \end{split}
  \end{equation*}
  Here $\mu$ and $J$ are given by the constraint equations $16 \pi \mu
  = \Scal - |K|_g^2 + (\tr K)^2$, and $8\pi J = \nabla^M\tr K - \div^M
  K$ and $(K^T)^\circ$ denotes the trace free part of the tangential
  projection of $K$ onto $\Sigma$, i.e. $(K^T)^\circ_{\alpha\beta} =
  K_{\alpha\beta} - \frac{1}{2}
  \gamma^{\eps\delta}K_{\eps\delta}\gamma_{\alpha\beta}$. Moreover,
  $\theta = K(\cdot,\nu)^T$, and $G$ denotes the Gaussian curvature of
  $\Sigma$.
\end{proposition}
\begin{myproof}
  Multiply \eqref{eq:lin1} with $f$ and integrate to obtain
  \begin{multline*}
    \int_\Sigma f \LHP f\dmu \\
    = \int_\Sigma |\nabla f|^2 - f^2 \big(|A|^2
    +\Ric(\nu,\nu) \mp \nabla_\nu \tr K 
    \pm \nabla_\nu   K(\nu,\nu)\big) 
    \pm 2f K(\nabla^\Sigma\! f,\nu)\dmu\,.
  \end{multline*}
  By the Gauss equation and the constraint equation we compute
  \begin{equation*}
    \label{eq:lin3}
    |A|^2 + \Ric(\nu,\nu) = 8\pi\mu+ \frac{1}{2}\left(|K|^2 - (\tr K)^2 +
    H^2 + |A|^2\right) - G\,.
  \end{equation*}
  Considering the term $2\int_\Sigma f K(\nabla^\Sigma
  f,\nu)\,\rmd\mu$, we obtain by partial integration that
  \begin{multline*}
    \label{eq:lin4}
    2\int_\Sigma f K(\nabla^\Sigma f,\nu)\dmu
    \\
    = -\int_\Sigma f^2 (-8\pi J(\nu) + \nabla^M_\nu \tr K - \nabla^M_\nu K(\nu,\nu) -H K(\nu,\nu)
    + K^T\cdot A )\dmu \,.
    \end{multline*}
  This gives the asserted identity in view of $|K|^2=|K^T|^2+2|\theta|^2+K(\nu,\nu)^2$.
\end{myproof}
This expression can be used to prove positivity. In the sequel we will
restrict ourselves to data $(M,g,K)$ which are
$(m,0,\sigma,\eta)$-asymptotically flat. By eventually increasing
$\sigma$, every set of $(m,\delta,\sigma,\bar\eta)$-asymptotically flat data can
be made $(m,0,\sigma,\eta)$-asymptotically flat for any choice of
$\eta>0$.
\begin{proposition}
  \label{prop:positiv}
  For $m>0$ and constants $C_1^A$, $C_2^A$, and $C_3^A$, there are
  $\eta_0=\eta_0(m,C^A)$ and $r_0=r_0(m, \sigma,
  C^A)$ such that if the data $(g,K)$ are
  $(m,0,\sigma,\eta_0)$-asymptotically flat and $\Sigma$ satisfies
  \eqref{eq:Heq}, conditions \eqref{eq:B1}--\eqref{eq:B4} as well as
  $\rmin>r_0$, then there is $\mu_1$ with
  \[ \mu_1 \geq 6m R_e^{-3} - o(1) R_e^{-3} + C(m,C^A) R_e^{-3-2/3} \,,\]
  such that for all functions $f\in C^\infty(\Sigma)$ with
  $\int_\Sigma f \dmu=0$ the following inequality holds
  \[ \mu_1 \int_\Sigma f^2\,\rmd\mu \leq \int_\Sigma f
  \LHP f\,\rmd\mu \,.\] 
  Here $o(1)$ is as described at the beginning of section \ref{s:apriori2}.
\end{proposition}
\begin{myproof}
  It is a well-known fact that a lower bound on the Gauss curvature
  $G\geq \kappa$ of a surface gives a lower bound $\lambda_1\geq 2
  \kappa$ on the first eigenvalue of its Laplace-Beltrami
  operator. This bound is provided by theorem \ref{thm:apriori}, such
  that for all $f$ with $\int_\Sigma f\dmu = 0$ we obtain
  \[ \left(\tfrac{1}{2}\bar H^2 + \tfrac{4m}{R_e^3} -
  o(1) R_e^{-3} + C(m,C^A) R_e^{-3-2/3} \right) \int_\Sigma f^2\dmu
  \leq \int_\Sigma |\nabla f|^2\dmu\,. \]
  From proposition \ref{prop:lin_pi}, the asymptotics of $K$, the
  $\sup$-estimates for $\Acirc$ from proposition \ref{prop:rund_sup},
  and the expression for $G$ in theorem \ref{thm:apriori}, we obtain
  \begin{equation*}
    \begin{split}
    \int_\Sigma f \LHP f\,\rmd\mu 
    &\geq 
    \int_\Sigma (|\nabla f|^2 - |\Acirc|^2 
    - \tfrac{3}{4} H^2 + G - o(1) R_e^{-3} ) f^2 \dmu 
    \\
    &\geq
    \left( 6m R_e^{-3} - o(1)R_e^{-3} -
    C(m,C^A)R_e^{-3-2/3}\right)\int_\Sigma f^2\dmu\,.
    \end{split}
  \end{equation*}
  If $o(1)$ is small enough, the factor  on the right hand side is
  positive, and this gives the assertion.
\end{myproof}
We are now able to show that solutions $u$ of $\LHP u=\const$ are
almost constant.
\begin{proposition}
  \label{prop:lin_sup}
  Let $(M,g,K)$ and $\Sigma$ be as in proposition
  \ref{prop:positiv}. Consider a solution $u$ of $\LHP u = f$ with
  $ \int_\Sigma(f-\bar f)^2\dmu \leq \mu_1^2/4\bar u^2 $ where $\mu_1$
  is as in proposition \ref{prop:positiv}, $\bar f=
  |\Sigma|^{-1}\int_\Sigma f\dmu$ is the mean value of $f$, and $\bar
  u$ is the mean value of $u$. Then
  \[ \sup_\Sigma |u-\bar u| \leq \left(o(1) + C(m,C^A)R_e^{-2/3}\right) \bar{u}\,. \]
\end{proposition}
\begin{myproof}
  We can assume that $u$ is normalized such that $\bar u =1$. Then
  \[ \LHP (u-1) = f + \left(|A|^2 +\Ric(\nu,\nu) \pm
  \nabla^M_\nu K(\nu,\nu) \mp \nabla^M_\nu\tr K \right)\,. \]
  Multiplying by $(u-1)$, integrating, and using proposition
  \ref{prop:positiv}, we obtain
  \begin{multline*}
    \mu_1 \int_\Sigma (u-1)^2 \dmu
    \\
    \leq
    \int_\Sigma (u-1)f + (u-1)\left(|A|^2 +\Ric(\nu,\nu) \pm
    \nabla^M_\nu K(\nu,\nu) \mp \nabla^M_\nu\tr K \right)\dmu\,.
  \end{multline*}
  Using the Schwarz inequality and the assumption on $f$ we estimate
  \begin{equation}
    \label{eq:lin_sup1}
    \int_\Sigma (u-1) f\dmu 
    = \int_\Sigma (u-1) (f-\bar f)\dmu 
    \leq \frac{\mu_1}{2} \left(\int_\Sigma (u-1)^2\dmu\right)^{1/2}\,.
  \end{equation}
  Define $R_e$ and $\bar H$ as in theorem \ref{thm:apriori}, then
  \begin{equation*}
    |\Acirc|^2 + \tfrac{1}{2}|H^2-\bar H^2| + |\Ric(\nu,\nu) +
    2mR_e^{-3}|
    \leq o(1)R_e^{-3} + C(m,C^A)R_e^{-3-2/3}\,.
  \end{equation*}
  Combining $\int_\Sigma(u-1)\bar H^2\dmu=0$ with
  the Schwarz inequality gives 
  \begin{multline*}
    \left|\int_\Sigma (u-1) \left(|A|^2 +\Ric(\nu,\nu) \pm
      \nabla^M_\nu K(\nu,\nu) \mp \nabla^M_\nu\tr K
      \right)\,\rmd\mu\right|
      \\
      \leq
      \left( o(1)R_e^{-2} + C(m,C^A) R_e^{-2-2/3}\right)
      \|u-1\|_{L^2(\Sigma)} \,.
  \end{multline*}
  Inserting this into \eqref{eq:lin_sup1}, we obtain the $L^2$-estimate
  \begin{equation*}
    \|u-1\|_{L^2(\Sigma)}^2 
    \leq \mu_1^{-2} \left( o(1)R_e^{-4} + C(m,C^g)R_e^{-4-4/3}\right) \,.
  \end{equation*}
  By standard estimates from the theory of linear elliptic partial
  differential equations of second order \cite{GT:1998} we can obtain a
  $\sup$-estimate from this
  \[ \sup_\Sigma |u-1| \leq \mu_1^{-1}\left(o(1)R_e^{-3}+ C(m,C^A)
  R_e^{-3-2/3}\right)\,, \] 
  which implies the assertion, in view of the estimate for $\mu_1$
  from proposition \ref{prop:positiv}.
\end{myproof}
\begin{corollary}
  \label{koro:vorzeichen}
  Provided $o(1)$ is small enough, and $f$ is as in the previous
  proposition, a solution of $Lu=f$ does not change sign.
\end{corollary}
\begin{corollary}
  \label{koro:lin_sup}
  Let $u$ be a solution of $\LHP u = f$. If
  \[ \int_\Sigma (u-\bar u)f\dmu \leq \frac{\mu_1}{2} \left(\int_\Sigma(u-\bar u)^2\dmu\right)^{1/2}\,, \]
  with $\mu_1$ from proposition \ref{prop:positiv}, then
  \[ \sup_\Sigma |u-\bar u| \leq o(1) + C(m,C^A)R_e^{-2/3} \bar{u}\,.\] 
\end{corollary}
This corollary implies that $\LHP$ is invertible in  suitable Banach
spaces.
\begin{theorem}
  \label{thm:linearisierung}
  Under the assumptions of the previous proposition, $L^{\CH\pm\CP}$
  is invertible as operator $L^{\CH\pm\CP} : C^{2,\alpha}(\Sigma)
  \rightarrow C^{0,\alpha}(\Sigma) $ for each $0<\alpha<1$. Its
  inverse $ \LHP_\text{inv} : C^{0,\alpha}(\Sigma) \rightarrow
  C^{2,\alpha}(\Sigma) $ exists and is continuous. It satisfies $\|
  \LHP_\text{inv} f \|_{L^2(\Sigma)} \leq R_e^3/3m\,
  \|f\|_{L^2(\Sigma)}$ and the H\"older norm estimate
  \[ \| \LHP_\text{inv} f \|_{C^{2,\alpha}(\Sigma)} \leq 
  C(\alpha,\Sigma) \tfrac{R_e^3}{3m}\|f\|_{C^{2,\alpha}(\Sigma)} \,.\]
\end{theorem}
\begin{myproof}
  Assume that there exists a function $u$ with $\|u\|_{L^2(\Sigma)} = 1$ 
  and
  \begin{equation}
    \label{eq:lin_kern}
    \sup_{\|v\|_{L^2(\Sigma)}=1}\ \left|\int_\Sigma v \LHP u \dmu \right|\leq
    \frac{3m}{R_e^3}\,.
  \end{equation}
  From proposition \ref{prop:positiv} we have that $\bar u \neq
  0$. Without loss of generality, $\bar u>0$.  Choosing  $v=u-\bar
  u$ in \eqref{eq:lin_kern} implies that the assumptions of corollary
  \ref{koro:lin_sup} are satisfied. If $o(1)$ is small enough, we
  obtain that $\bar u/2 \leq u \leq 2\bar u$. From $\| u
  \|_{L^2(\Sigma)}=1$ we obtain that $\bar u\geq
  \frac{1}{2}|\Sigma|^{-1/2}$, and from H\"older's inequality $\bar
  u \leq |\Sigma|^{-1/2}$. Using $v=1$ in \eqref{eq:lin_kern} gives 
  \begin{equation}
    \label{eq:lin_test1}
    \left|\int_\Sigma \LHP u \dmu\right| \leq \frac{3m}{R_e^3}|\Sigma| \leq
    C(m,C^A) R_e^{-1}\,.
  \end{equation}
  On the other hand, we compute from \eqref{eq:lin1}, using partial
  integration, that
  \begin{equation*}
    \begin{split}
      \label{eq:lin_eq}
      \int_\Sigma \LHP u\dmu &= - \int_\Sigma u \big(|A|^2 +
      \Ric(\nu,\nu) \pm \nabla_\nu^MK(\nu,\nu) \mp \nabla_\nu^M\tr K
      \\
      &\qquad\qquad \pm \nabla^M_{e_\alpha}K(e_\alpha,\nu)\mp HK(\nu,\nu)\pm
      K^T\cdot A\big)\dmu \,.
    \end{split}
  \end{equation*}
  Inserting this into the previous estimate, we infer using \eqref{eq:lin_test1}, that
  \begin{equation*}
    \left|\int_\Sigma u|A|^2 \dmu\right| \leq \left|\int_\Sigma
    \LHP u\dmu\right| + C(m)R_e^{-3}|\Sigma|\bar u \leq C(m)R_e^{-1}\bar u
    \,.
  \end{equation*}
  From $\bar u \leq 2 u$ we obtain that $\int_\Sigma H^2\dmu \leq C(m)R_e^{-1}$,
  which contradicts \eqref{eq:B2} for large $R_e$. This implies that
  $\LHP$ is injective, and since it is a linear elliptic
  operator, the Fredholm alternative consequently implies its
  surjectivity. The existence of a continuous inverse $\LHP_\text{inv}$
  with the asserted bounds follows \cite[Chapter 5]{GT:1998}. Note
  that by the a priori estimates of theorem \ref{thm:apriori} the
  Gauss curvature, and therefore the injectivity radius, are controlled.
\end{myproof}
\begin{remark}
  The constant $C(\alpha,\Sigma)$ can be chosen uniformly by using 
  Schauder estimates in e.g. harmonic coordinate patches on
  $\Sigma$.
  Analogous estimates in the spaces $W^{2,p}(\Sigma)$ can be found in
  \cite[Chapter 2]{CK:1993}
  \[ \| \LHP_\text{inv} f \|_{W^{2,p}(\Sigma )} \leq 
  C(2,p) \tfrac{R_e^3}{3m}\|f\|_{L^p(\Sigma )} \,.\]
  The constants $C(2,p)$ therein can be chosen uniformly since they only depend on
  $k_\text{min}:=|\Sigma|^{-1}\min_\Sigma G$, and
  $k_\text{max}:=|\Sigma|^{-1}\max_\Sigma G$, which are controlled in
  our case.
\end{remark}

  \section{The foliation}
\label{s:foliation}
To prove the existence of surfaces satisfying $H\pm P =
\const$, we use the following strategy. Let $(g,K)$ be
$(m,0,\sigma,\eta)$-asymptotically flat with $m>0$. Let
$g_\tau:=(1-\tau)g^S+\tau g$, and $K_\tau:= \tau K$. Then the data
$(g_\tau,K_\tau)$ is also $(m,0,\sigma,\eta)$-asymptotically
flat. For the initial reference data $(g^S,0)$ we know a lot of
solutions to the equation $H=\const$, namely the centered spheres
(note that if $K\equiv 0$ then $P\equiv 0$). The mean curvature of a
centered sphere of radius $r$ in with respect to $g^S$ can be computed using
\ref{lemma:extcurv} and equals
\[ H^S(r) = \left(1+\frac{m}{2r}\right)^{-3}\left(\frac{2}{r} -
\frac{m}{r^2}\right)\,. \] This function is invertible for
$r>r_1(m)$. The inverse function satisfies $|r-2/h|\leq C$, for any $C$
provided $h_1$ is chosen small enough. Let $h>0$ be a constant. Then
we can solve $H^S(r)=h$ with $r>r_1(m)$, provided $h<h_1(m)$. Therefore
the equation
\[ H \pm P = h \]
is satisfied on a sphere of radius $r(h)$ for $\tau = 0$.  To deform
this solution for $\tau=0$ to a family of solutions for
$\tau\in[0,1]$, we introduce two classes of surfaces. For this
consider the following conditions related to \eqref{eq:B1}--\eqref{eq:B4} by
appropriately choosing the constants
\begin{gather}
  \label{eq:C1} \tag{B1}
  R(\Sigma) \leq 8 \rmin\,, \\
  \label{eq:C2} \tag{B2}
  R(\Sigma)^{-1} \leq 8(H \pm P)\,, \\
  \label{eq:C3} \tag{B3}
  \int_\Sigma u\,|A|^2 \dmu \leq 8 \int_\Sigma u\,\det A\dmu \qquad\text{for
    all}\qquad 0 \leq u\in C^\infty(\Sigma)\,, \\
  \label{eq:C4} \tag{B4}
  |\Sigma|_e^{-1} \int_\Sigma \id_\Sigma \dmu^e \leq \frac{3}{4} R_e
\end{gather}
Choose $\eta_0$ so small, and $r_0$ so large, that corollaries
\ref{koro:B1'}, \ref{koro:B2'}, and theorem \ref{thm:apriori}
imply that these conditions hold with better constants on surfaces
$\Sigma$ with $\rmin>r_0$
\begin{gather}
  \label{eq:D1} \tag{C1}
  R(\Sigma) \leq 4 \rmin \,,\\
  \label{eq:D2} \tag{C2}
  R(\Sigma)^{-1} \leq 4(H \pm P)  \,,\\
  \label{eq:D3} \tag{C3}
  \int_\Sigma u\,|A|^2 \dmu \leq 4 \int_\Sigma u\,\det A\dmu \qquad\text{for
  all}\qquad 0 \leq u\in C^\infty(\Sigma)\,. \\
  \label{eq:D4} \tag{C4}
   |\Sigma|_e^{-1} \int_\Sigma \id_\Sigma \dmu^e \leq \frac{7}{8} R_e
\end{gather}
By eventually decreasing $\eta_0$ and increasing $r_0$, we can assume
that \eqref{eq:D1} -- \eqref{eq:D4} imply that the linearized operator
$\LHP$ from the previous section is invertible, corollary
\ref{koro:gradient} guarantees that $\Sigma$ is globally a graph over
$S^2$, and $g^e(\nu^e,\rho)>1/2$. Moreover, from theorem
\ref{thm:apriori} we can assume that for all surfaces satisfying
\eqref{eq:C1} -- \eqref{eq:C4}, also
\begin{equation}
  \label{eq:ex_radius}
  \frac{1}{4} \rmin \leq (H\pm P)^{-1} \leq 4\,\rmin \,.
\end{equation}
Let $(g,K)$ be data such that for fixed $m>0$ the data
$(g_\tau,K_\tau)$ as before all are
$(m,0,\sigma,\eta_0)$-asymptotically flat. Define the following
nested sets of surfaces:
\begin{equation*}
  \begin{split}
  \CS_1(\tau) &:= \{ S^2\approx\Sigma \subset M :\Sigma\ \text{satisfies}\ \rmin > r_0\ \text{and
  \eqref{eq:C1}--\eqref{eq:C4} w.r.t.}\,(g_\tau,K_\tau) \}\\
  \CS_2(\tau) &:= \{ S^2\approx\Sigma \subset M :\Sigma\ \text{satisfies}\ \rmin > 2r_0\ \text{and
  \eqref{eq:D1}--\eqref{eq:D4} w.r.t.}\,(g_\tau,K_\tau) \}
  \end{split}
\end{equation*}
Choose $0<h_2\leq h_1$ such that the centered spheres $S_r(0)$ with
mean curvature $H<h_2$ are in $\CS_2(0)$. Choose
$h_0<\min\{h_1,h_2,\frac{1}{8}r_0^{-1}\}$. Let
\[ \kappa:[0,1]\to(0,h_0) \times [0,1] : t \mapsto \left(h(t),
\tau(t)\right)\] 
be a continuous, piecewise smooth curve with $\tau(0)=0$. Denote by
$(H\pm P)_\tau$ the function $H\pm P$ evaluated with respect to
$(g_\tau,K_\tau)$. Let $I_\kappa \subset[0,1]$ be the set
\[ I_\kappa:=\left\{ t\in[0,1]: \exists\,\Sigma(t)\in\CS_2(\tau(t))\ \text{with}\ (H\pm P)_{\tau(t)}=h(t)\right\} \]  
\begin{proposition}
  Under the assumptions of this section, $I_\kappa=[0,1]$. 
\end{proposition}
\begin{myproof}
  We can assume that $\kappa$ is smooth. By choice of $h_0$,  $0\in
  I_\kappa$, so $I_\kappa$ is nonempty. 
 
  For proving that $I_\kappa$ is open, let $t_0\in I_\kappa$, and
  $\Sigma\in \CS_2(\tau(t_0))$ the surface with $(H\pm
  P)_{\tau(t_0)}=h(t_0)$. Consider Gaussian normal coordinates $y:
  \Sigma\times (-\eps,\eps)\to M$, and let $B:=\{f\in C^{2,\alpha}(\Sigma):
  \sup |f|<\eps\}$. Define the operator
  \[ \CL : B \times [0,1] \to C^{0,\alpha}(\Sigma) : (f,t) \mapsto (\CH\pm\CP)_t(f) - h(t)\,, \]
  where $(\CH\pm\CP)_t(f)$ is the function $(H\pm P)_t$ on
  $\graph(f)$. This operator is differentiable, and we have $\CL(0,t_0)
  = 0$.
  
  The differential of $\CL$ with respect to the first variable is the
  operator $\LHP$ from section \ref{s:linearized}, and is invertible since
  $\Sigma\in\CS_2(\tau(t))$. By the implicit function theorem there
  exists $\delta>0$, and a differentiable function
  $\xi:(t_0-\delta,t_0+\delta)\to B$, such that $\CL(\xi(t),t)=0$ for
  all $t$ with $|t-t_0|<\delta$.
  
  Hence, for each such $t$ there is a surface $\Sigma(t)$ with $(H\pm
  P)_{\tau(t)}=\const$. By continuity, and by eventually decreasing
  $\delta$, we can assume that $\Sigma(t)\in\CS_1(\tau(t))$. By choice
  of $r_0$ and $\eta_0$ conditions \eqref{eq:C1}--\eqref{eq:C4} imply
  \eqref{eq:D1}--\eqref{eq:D4}. By choice of $h_0$ we obtain
  $\rmin>2r_0$ whence $\Sigma(t)\in \CS_2(\tau(t))$. That is,
  $I_\kappa$ contains a small neighborhood of $t_0$.
  
  To show that $I_k$ is closed, assume that $\{t_n\}\subset I_\kappa$
  is a convergent series with $\lim_{n\to\infty}t_n\to t$. Let
  $\Sigma(t_n)\in \CS_2(\tau(t_n))$ be the surface with $(H\pm
  P)_{\tau(t_n)}=h(t_n)$. By corollary \ref{koro:gradient} all 
  $\Sigma(t_n)=\graph(u_n)$ are graphs over $S^2$ as described in
  section~\ref{s:linearized}.
  
  From the position estimates in proposition \ref{prop:pos}, the uniform
  estimates for the angle $g^e(\nu^e,\rho)$, and the uniform curvature estimates
  from corollary \ref{koro:B2'} and proposition \ref{prop:rund_sup} we
  obtain uniform $C^2(S^2)$-estimates for the sequence $(u_n)$. In addition, the
  $W^{1,2}$-estimates on the curvature imply uniform $W^{3,2}$-estimates
  for $(u_n)$.
  
  We can assume that $(u_n)$ converges in $W^{2,p}(S^2)$ to $u\in
  W^{2,p}(S^2)$ for a $1<p<\infty$. Furthermore, we can assume that
  $(u_n)\to u$ in $C^{1,\alpha}(S^2)$ for a fixed $0<\alpha<1$.
  
  On $\graph(u)$ a weak version of the quasilinear equation $(H\pm
  P)_{\tau(t)} = h(t)$ is satisfied. By fixing coefficients, this can be
  interpreted as a linear equation. Since $u\in C^{1,\alpha}$, the
  coefficients of this equation are $C^{0,\alpha}$. Regularity
  theory for such equations \cite[Chapter 8]{GT:1998} implies that
  $u$, and therefore $\Sigma$, are smooth. By $C^{1,\alpha}$-convergence
  $\Sigma$ satisfies \eqref{eq:D1}, \eqref{eq:D4}, and $\rmin>2 r_0$. By
  $W^{2,p}$-convergence \eqref{eq:D2} and \eqref{eq:D3} are satisfied,
  provided $p$ is large enough. Therefore $t\in I_\kappa$, and
  $I_\kappa$ is closed.
\end{myproof}
This gives the following:
\begin{theorem}
  \label{thm:existenz}
  Let $m>0$ be fixed. Then there exists $h_0=h_0(m,\sigma)$ and
  $\eta_0=\eta_0(m)$ such that for every
  $(m,0,\sigma,\eta_0)$-asymptotically flat data set $(g,K)$ and every
  curve $\kappa:[0,1]\to (0,h_0)\times[0,1]:t\mapsto (h(t),\tau(t))$
  there exists a smooth family of surfaces $\Sigma_\kappa(t)\in
  \CS_2(\tau(t))$ satisfying $H\pm P= h(t)$ with respect to the
  $\tau(t)$-data.
\end{theorem}
\begin{remark}
  \label{rem:eindeutigkeit}
  At first glance, the resulting $H\pm P=\const$-surface could depend on
  the choice of the curve $\kappa$ from $\kappa(0)$ to
  $\kappa(1)$. However, since the range of $\kappa$ is simply connected,
  and the solutions obtained from the implicit function theorem are
  locally unique, a standard argument using the homotopy of two curves
  with common endpoints shows that the surfaces in fact only depend on
  the endpoints of $\kappa$.
\end{remark}
We are now ready to prove the existence part of theorem
\ref{thm:haupt}.
\begin{theorem}
  \label{thm:haupt_exist}
  Let $m>0$ be fixed and $\eta_0$ and $h_0$ be as in theorem
  \ref{thm:existenz}. By possibly decreasing $\eta_0$ and $h_0$ we
  assure that corollary \ref{koro:vorzeichen} is valid. Then the
  surfaces satisfying $H\pm P =\const$ constructed in theorem
  \ref{thm:existenz} form a foliation. For small $H\pm P$ these
  surfaces have arbitrarily large radius. In addition, there is a
  differentiable map
  \[\CF:S^2\times(0,h_0)\times [0,1] \to M\]
  such that the surfaces $\CF(S^2,h,\tau)$ satisfy $H\pm P = h$ with
  respect to the data $(g_\tau,K_\tau)$. This foliation can therefore be
  obtained by deforming a piece of the $H=\const$ foliation of
  $(\IR^3,g^S)$ by centered spheres.
\end{theorem}
\begin{myproof}
  Choose $0<h<h_0$, and define the curve $\kappa_h(t)=(h,t)$ for
  $t\in[0,1]$. Using theorem \ref{thm:existenz} we obtain a family of surfaces
  $\Sigma_{h,\tau}$ with $H\pm P=h$ by deforming the centered sphere which
  has $H^S=h$ with respect to $g^S$ along $\kappa$. The position
  estimates and \eqref{eq:ex_radius} imply $h^{-1} \leq 4
  \rmin(\Sigma_{h,t})$, such that by choosing $h$ small, we can make
  $\rmin$ of $\Sigma_1(h)$ large.
  
  The map $\CF$ can be constructed by setting
  $\CF(S^2,h,\tau)=\Sigma_{h,\tau}$ and defining the parametrization
  of $\Sigma_{h,\tau}$ by the fact that $\Sigma_{h,\tau}$ is a graph
  over $S^2$. This implies the differentiability of $\CF$ with respect
  to $p\in S^2$ and $\tau\in[0,1]$.

  To show that these surfaces form a foliation, choose another
  curve. Let $h_1\in(0,h_0)$ be fixed. The curve $\kappa_{h_1}$ gives
  a fixed reference surface $\Sigma_{h_1,1}$. For $h_2<h_1$ consider the
  curves $\lambda_{h_2}(t) = ((1-t)h_1+ th_2,1)$. Concatenating
  $\kappa_{h_1}$ and $\lambda_{h_2}$ gives a family of surfaces
  $\Sigma'_{h,1}$ with $h\in[h_2,h_1]$ as well as a differentiable
  map $F : S^2 \times [h_2,h_1] \to M$ such that
  $F(S^2,h)=\Sigma'_{h,1}$. Remark \ref{rem:eindeutigkeit} implies
  that $\Sigma_{h,1}'=\Sigma_{h,1}=:\Sigma_h$. Therefore $\CF$ is
  differentiable with respect to $h\in(0,h_0)$.

  Let $\nu_h$ denote the normal to $\Sigma_{h}$, then the lapse
  $\alpha_h$ of the family $F$ is defined as $\alpha_h:=
  g\left(\nu_h,\DD{F}{h}\right)$. Since $H\pm P=\const$ along
  $\Sigma_h$, and therefore the tangential part of $\DD{F}{h}$ is
  irrelevant for the evolution of $H\pm P$, we have
  \[ h_1-h_2 = \DD{}{h}(H\pm P) = \LHP\alpha_h \]
  with the operator $\LHP$ from section \ref{s:linearized}. By
  corollary \ref{koro:vorzeichen}, $\alpha_h$ does not change
  sign. Therefore the family of the $\Sigma_h$ is a foliation.
\end{myproof}
We can also prove the uniqueness of $H\pm P=\const$ surfaces.
\begin{theorem}
  For $m>0$ there are $\eta_0(m,C^A)>0$ and $h_0(m,\sigma,C^A)>0$ such
  that if $(g,K)$ is $(m,0,\sigma,\eta_0)$-asymptotically flat, then
  two surfaces $\Sigma_1$ and $\Sigma_2$ satisfying
  \eqref{eq:B1}--\eqref{eq:B4} and $H\pm P = h =\const$ with
  $h\in(0,h_0)$ coincide.
\end{theorem}
\begin{myproof}
  We prove this by reversing the process we used in the proof of the
  existence result. That is, we start for the data $(g_\tau,K_\tau)$
  at $\tau=1$ with $\Sigma_1$ and $\Sigma_2$ and obtain surfaces
  $\Sigma'_1$ and $\Sigma'_2$ with $H=h=\const$ with respect to the
  Schwarzschild metric $g^S$ at $\tau=0$. Here $\eta_0$ and $h_0$ have
  to be adjusted as in the beginning of this section, such that this
  process works.
  
  By the uniqueness of such surfaces satisfying
  \eqref{eq:B1}--\eqref{eq:B4} in the Schwarzschild metric, as follows
  for example from Huisken and Yau \cite[Section 5]{HY:1996}, we infer
  that $\Sigma'_1$ coincides with $\Sigma'_2$. Then by the local
  uniqueness of the implicit function theorem also $\Sigma_1$ and
  $\Sigma_2$ coincide.
\end{myproof}
\begin{corollary}
  The $H\pm P =\const$ foliations from theorem \ref{thm:haupt}
  consisting of surfaces satisfying \eqref{eq:B1}--\eqref{eq:B4} are
  unique at infinity.
\end{corollary}

  \section{Special data}
\label{s:special}
We want to interpret the foliation of $H\pm P=\const$ surfaces in a
physical manner. A foliation of surfaces satisfying $H=\const$ was
interpreted in \cite{HY:1996} as the center of mass of an isolated
system. The definition of this foliation does not refer to the
extrinsic curvature $K$ and therefore can not contain information on
dynamical physics. In contrast, proposition \ref{prop:york_lin} shows that
the $H\pm P=\const$ foliation allows an interpretation as linear
momentum.

We restrict ourselves to data $(g,K)$ with
$\|g-g^S\|_{C^2_{-1-\delta}}<\infty$ with $\delta>0$ and 
\[ K = \frac{3}{2r^2}\big(\rho\tensor p + p\tensor \rho - 2\la
p,\rho\ra(g^e- \rho\tensor\rho)\big) + O(r^{-2-\delta})\, \] where
$p\in\IR^3$ is a fixed vector, $\rho=x/r$ is the radial direction, and
the derivatives of $O(r^{-2-\delta})$ are of order
$O(r^{-3-\delta})$. This structure of $K$ was proposed by York
\cite{York:1978} and represents a trace free extrinsic curvature
tensor with ADM-momentum $p$. There exist initial data satisfying the
constraint equations with these asymptotics.
Using this representation of $K$, we can refine the 
estimates from proposition \ref{prop:pos} and obtain
\begin{proposition}
  \label{prop:york_lin}
  Let $(g,K)$ be as described above. If $|p|<m$ is small enough, and
  $\Sigma$ satisfies $H\pm P=\const$, assumptions
  \eqref{eq:B1}--\eqref{eq:B4}, and $\rmin>r_0$, then there exist a
  vector $a\in\IR^3$, a sphere $S=S_{R_e}(a)$, and a
  parameterization $\psi:S\to\Sigma$ such that
  \begin{eqnarray*}
    \left| a/R_e \mp \tau(v) \bar p \right| &\leq& CR_e^{-\delta} \,,\\
    \sup_S | \phi - \id_S | &\leq& CR_e^{-\delta} \,,\\
    \sup_\Sigma | \nu^e - R_e^{-1} (r\rho - a) | &\leq& CR_e^{-\delta}
  \end{eqnarray*}
  with $\bar p = \frac{p}{|p|}$, $v=\frac{|p|}{m}$, and $\tau(v)=\frac{1-\sqrt{1-v^2}}{v}$. If
  $0\leq v\leq 1$ then $0\leq \tau(v) \leq 1$ and
  $\tau(v) = \frac{1}{2}v + O(v^3)$ for $v\to 0$.
\end{proposition}
\begin{myproof}
  This proof is similar to the proof of proposition
  \ref{prop:pos}. However, instead of estimating like
  \eqref{eq:pos_p_est}, we compute more carefully using the asymptotics
  of $K$. For the test vector $b=\bar a := \frac{a}{|a|}$ we obtain
  \begin{equation}
    \label{eq:postest_a}
    \left| - 8\pi m \frac{|a|}{R_e} \mp 4\pi \left\la p,\frac{a}{|a|}\right\ra \frac{|a|^2 +
      R_e^2}{R_e^2} \right| \leq CR_e^{-\delta}\,.
  \end{equation}
  Now we split $p=g^e(\bar a, p) \bar a + g^e(\bar q, p)\bar q$ with $\bar q$
  orthogonal to $a$ and $|\bar q|=1$. Then we use $\bar q$ as an additional test
  vector. This gives the second estimate
  \begin{equation}
    \label{eq:postest_q}
    |\la p, \bar q\ra| \left|\frac{4\pi}{5}\frac{5R_e^3 - 2|a|^2R_e - 3|a|^2}{R_e^3}\right|\leq
    CR_e^{-\delta}\,.
  \end{equation}
  Proposition \ref{prop:pos} gives $\tau:=|a|/R_e<1$ if $p$ is
  small. Then \eqref{eq:postest_q} implies that $|g^e(p,\bar q)| \leq
  CR_e^{-\delta}$, and therefore $|g^e(p,\bar a)-|p||=|g^e(p,\bar
  q)|\leq CR_e^{-\delta}$. Using \eqref{eq:postest_a} we infer that
  \[ \left| - 2m\tau \mp |p|(1+\tau^2) \right| \leq CR^{-\delta}\,, \]
  which implies the proposition.
\end{myproof}
\begin{remark}
  \begin{mynum}
  \item This means that surfaces $\Sigma(h)$ satisfying $H\pm
    P=h=\const$ are not only increasing in size for $h\to 0$, but that
    they also translate. The magnitude of this translation can be used to compute
    $p$. The asymptotic translation $\tau$ from the previous
    proposition can be found by comparing the Euclidean center of
    gravity to the center of gravity computed using the
    $g$-metric. In particular \[ \lim_{h\to 0} \left(|\Sigma^e|^{-1}
    \int_{\Sigma(h)} x \dmu^e - |\Sigma|^{-1} \int_{\Sigma(h)} x
    \dmu\right) = \frac{2}{3} m \tau \bar p\,. \]
    Here $\tau=\lim_{h\to 0} |a|/R_e(\Sigma(h))$ is the limit of the
    magnitude of the translation vector and $\bar p$ the unit vector
    pointing into its direction. Then $p$ can be
    computed from
    \[ \pm p = \frac{2m\tau}{1+\tau^2}\,\bar p \,.\]
  \item Corvino and Schoen \cite{CS:2003} also propose a standard form
    of the extrinsic curvature tensor, namely
    \[ K^\text{CS} = \frac{2}{r^2} \left( p \tensor \rho + \rho \tensor p - \la p, \rho \ra
    g^e \right) + O(r^{-3})\,.\] 
    Contrary to the York-form this is not
    trace free in the terms of highest order. Corvino and Schoen prove
    that data satisfying this asymptotic condition for $K$ and $g= g^S
    + O(r^{-2})$ are dense with respect to suitable, weighted Sobolev
    norms in the set of data $(\bar g, \bar K)$ satisfying the
    constraint equations and
    \[ \bar g = g^e + O(r^{-1}) \qquad\text{and}\qquad \bar K =
    O(r^{-2})\]
    Therefore these asymptotics posses a certain universality. 

    For these asymptotics we can also compute the asymptotic
    translation. It satisfies $\tilde\tau(v)=
    \frac{1-\sqrt{1-\frac{16}{15}v^2}}{\frac{8}{5}v}$. Here
    $\tilde\tau(v)=\frac{1}{3}v+ O(v^3)$ for $v\to 0$.

    This is not satisfactory for two reasons. At first, this
    asymptotic translation and the associated linear momentum formula
    do not coincide with the formula obtained from the York
    asymptotics. Secondly $0\leq \tilde\tau(v)\leq 1$ only for
    $0<v<\frac{15}{16}$, while from physical reasons at least the
    interval $v\in[0,1]$ should be admitted.
    
    On the other hand, we can not expect to obtain a valid formula
    independent of the slicing condition. For the $H\pm P=\const$
    foliation, therefore the slicing condition $\tr K = 0$ seems to be
    appropriate.
  \item
    Both the asymptotics of York and the asymptotics of Corvino and
    Schoen allow examples of initial data satisfying the vacuum
    constraint equations. This implies that the Sobolev norm used by
    Corvino and Schoen to prove density of data with their asymptotics
    is strong enough to preserve ADM-mass and -momentum, but not
    strong enough to reproduce the fine structure of the $H\pm
    P=\const$ foliation.\fillbox
  \end{mynum}
\end{remark}

  \bibliographystyle{amsalpha}
  \bibliography{references}
\end{document}